\documentclass[reqno,letterpaper]{amsart}

\usepackage{amsmath}
\usepackage{amsthm}

\usepackage{times}
\usepackage{bm}
\usepackage{natbib}
\usepackage{mathrsfs}
\usepackage{graphicx}
\usepackage[colorlinks=true,citecolor=blue]{hyperref}

\usepackage[plain,noend]{algorithm2e}

\makeatletter
\renewcommand{\algocf@captiontext}[2]{#1\algocf@typo. \AlCapFnt{}#2} %
\def\@algocf@capt@plain{top}
\renewcommand{\algocf@makecaption}[2]{%
  \addtolength{\hsize}{\algomargin}%
  \sbox\@tempboxa{\algocf@captiontext{#1}{#2}}%
  \ifdim\wd\@tempboxa >\hsize%
    \hskip .5\algomargin%
    \parbox[t]{\hsize}{\algocf@captiontext{#1}{#2}}%
  \else%
    \global\@minipagefalse%
    \hbox to\hsize{\box\@tempboxa}%
  \fi%
  \addtolength{\hsize}{-\algomargin}%
}
\makeatother

\usepackage[capitalize]{cleveref}  

\theoremstyle{plain}
\newtheorem{theorem}{Theorem}[section]
\newtheorem{proposition}[theorem]{Proposition}
\newtheorem{lemma}[theorem]{Lemma}
\newtheorem{corollary}[theorem]{Corollary}

\theoremstyle{definition}
\newtheorem{definition}[theorem]{Definition}

\theoremstyle{remark}
\newtheorem{remark}[theorem]{Remark}

\crefname{lemma}{Lemma}{Lemmas}
\crefname{proposition}{Proposition}{Propositions}
\crefname{corollary}{Corollary}{Corollaries}
\crefname{theorem}{Theorem}{Theorems}
\crefname{remark}{Remark}{Remarks}

\RequirePackage[usenames,dvipsnames]{color}
\usepackage{jhh-misc2}
\usepackage{autonum} 
\usepackage[no-extra-Laplace]{optional}

\newcommand{\Xvec}{X}

\newcommand{\targetdist}{\Pi}
\newcommand{\targetdensity}{\pi}
\newcommand{\approxdist}{\hat\targetdist}
\newcommand{\estdist}{\hat\targetdist}

\newcommand{\mainmeas}{\eta}
\newcommand{\estmeas}{\hat{\mainmeas}}
\newcommand{\mainPE}{U}
\newcommand{\estPE}{\hU}

\newcommand{\priordist}{\targetdist_{0}}
\newcommand{\postdist}{\targetdist}
\newcommand{\priordensity}{\targetdensity_{0}}
\newcommand{\postdensity}{\targetdensity}
\newcommand{\approxdensity}{\tilde\targetdensity}
\newcommand{\vardist}{Q}
\newcommand{\varfamily}{\mcQ}
\newcommand{\lik}{f}
\newcommand{\loglik}{\mcL}
\newcommand{\marginallik}{M}
\newcommand{\data}{Z}
\newcommand{\param}{\theta}
\newcommand{\Param}{\Theta}

\newcommand{\map}{\param^{\star}}
\newcommand{\maphessian}{H^{\star}}
\newcommand{\derivop}[1]{\grad^{#1}}
\newcommand{\hessian}{\derivop{2}}
\newcommand{\mapn}{\map_{n}}
\newcommand{\maphessiann}{H^{\star}_{n}}
\newcommand{\logpostn}{\ell_{n}}

\newcommand{\abscont}{\ll}  

\newcommand{\textder}[2]{\dee #1/\dee #2}

\newcommand{\Wiener}[1]{W_{#1}}
\newcommand{\dWt}[1]{\dee \Wiener{#1}}

\newcommand{\drift}{b}

\newcommand{\Xt}[1]{X_{#1}}
\newcommand{\dXt}[1]{\dee \Xt{#1}}
\newcommand{\Yt}[1]{Y_{#1}}
\newcommand{\dYt}[1]{\dee \Yt{#1}}
\newcommand{\dt}{\dee t}

\newcommand{\Lparg}[2]{L^{#1}(#2)}
\newcommand{\Lpnormarg}[3]{\staticnorm{#3}_{\Lparg{#1}{#2}}}

\newcommand{\couplings}[2]{\Gamma(#1, #2)}
\newcommand{\coupling}{\gamma}

\newcommand{\pwassSimple}[3]{\mcW_{#1}(#2, #3)}
\newcommand{\dchi}[2]{d_{\chi}(#1 || #2)}

\newcommand{\fd}[3]{d_{#1}(#2, #3)}

\newcommand{\MAD}[1]{\operatorname{MAD}_{#1}}

\newcommand{\mean}[1]{\mu_{#1}}
\newcommand{\meanvec}[1]{\mu_{#1}}

\newcommand{\var}[1]{\sigma_{#1}^{2}}
\newcommand{\std}[1]{\sigma_{#1}}

\def\norm#1{\left\|{#1}\right\|} 

\newcommand{\onenorm}[1]{\norm{#1}_1} 
\newcommand{\twonorm}[1]{\norm{#1}_2} 
\newcommand{\infnorm}[1]{\norm{#1}_{\infty}} 
\newcommand{\opnorm}[1]{\norm{#1}_{op}} 
\def\staticnorm#1{\|{#1}\|} 
\newcommand{\statictwonorm}[1]{\staticnorm{#1}_2} 
\newcommand{\inner}[2]{\langle{#1},{#2}\rangle} 
\newcommand{\binner}[2]{\left\langle{#1},{#2}\right\rangle} 

\def\Holder{H\"older\xspace}

\allowdisplaybreaks[3]

\begin{document}

\title[Practical bounds on posterior approximation error]{Practical bounds on the error of Bayesian posterior approximations: A nonasymptotic approach}

\author[J.~H.~Huggins]{Jonathan H.~Huggins}
\address{Department of Biostatistics, Harvard T.~H.~Chan School of Public Health}
\urladdr{http://jhhuggins.org}
\email{jhuggins@mit.edu}

\author[M.~Kasprzak]{Miko{\l}aj Kasprzak}
\address{Department of Statistics, University of Oxford}
\urladdr{https://mikolajkasprzak.wordpress.com/}
\email{kasprzak@stats.ox.ac.uk}

\author[T.~Campbell]{Trevor Campbell}
\address{Department of Statistics, University of British Columbia}
\urladdr{http://www.trevorcampbell.me/}
\email{trevor@stat.ubc.ca}

\author[T.~Broderick]{Tamara Broderick}
\address{Department of Electrical Engineering and Computer Science, MIT}
\urladdr{http://www.tamarabroderick.com}
\email{tbroderick@csail.mit.edu}

\begin{abstract}
Bayesian inference typically requires the computation of an approximation to the posterior distribution. 
An important requirement for an approximate Bayesian inference algorithm is to output high-accuracy posterior mean and uncertainty estimates. 
Classical Monte Carlo methods, particularly Markov Chain Monte Carlo, remain the gold standard for approximate Bayesian inference 
because they have a robust finite-sample theory and reliable convergence diagnostics.
However, alternative methods, which are more scalable or apply to problems where Markov Chain Monte Carlo cannot be used,
lack the same finite-data approximation theory and tools for evaluating their accuracy.  %
In this work, we develop a flexible new approach to bounding the error of mean and uncertainty estimates of scalable inference algorithms.
Our strategy is to control the estimation errors in terms of Wasserstein distance, then bound the Wasserstein distance via a generalized notion of Fisher distance.
Unlike computing the Wasserstein distance, which requires access to the normalized posterior distribution,
the Fisher distance is tractable to compute because it requires access only to the gradient of the log posterior density.
We demonstrate the usefulness of our Fisher distance approach by deriving bounds on the Wasserstein error 
of the Laplace approximation and Hilbert coresets. 
We anticipate that our approach will be applicable to many other approximate inference methods 
such as the integrated Laplace approximation, variational inference, and approximate Bayesian computation.

\end{abstract}

\maketitle

\section{Introduction}

Exact Bayesian statistical inference is known for providing point estimates with desirable decision-theoretic properties as well as coherent uncertainties. 
Using Bayesian methods in practice, though, typically requires approximations to the posterior distribution and hence to both point estimates and uncertainties.
Hence, it is crucial to quantify the error introduced by such approximations. 
Monte Carlo methods, especially Markov chain Monte Carlo, are the gold standard for approximate Bayesian inference in part due to their flexibility and 
strong theoretical guarantees on quality for finite data.
However, these guarantees are typically asymptotic in running time, and computational concerns have motivated a spate of alternative Bayesian approximations. 
These include the Laplace approximation~\citep{Schervish:1995} and the integrated nested Laplace approximation~\citep{Rue:2009,Rue:2017},
approximate Bayesian computation~\citep{Marjoram:2003,Marin:2011,Karabatsos:2018}, 
subsampling Markov chain Monte Carlo~\citep{Welling:2011,Korattikara:2014,Bardenet:2014,Alquier:2016,Teh:2016,Vollmer:2016}, 
consensus methods~\citep{Scott:2013,Rabinovich:2015,Srivastava:2015,Li:2017},
and variational approaches~\citep{Blei:2017} such as automatic differentiation variational inference~\citep{Ranganath:2014,Kucukelbir:2015}.
While these methods have empirically demonstrated computational gains on problems of interest, rigorous characterization of their finite-data 
approximation accuracy remains underdeveloped, though ongoing~\citep{Alquier:2016b,Alquier:2017,Ogden:2017,CheriefAbdellatif:2018,Wang:2018,Ogden:2018,Pati:2018}.
We aim to provide theoretical tools to help address this gap.
In particular, since practitioners often report point estimates in the form of the mean or median and uncertainties in the form of the variance, 
standard deviation, or quantiles~\citep{Robert:1994,Gelman:2013}, these are the quantities we focus on approximating well. 
A natural approach is to start by bounding a statistical divergence between the exact and approximate posterior distributions, 
then use this bound to in turn bound the error in approximate posterior functionals of interest.

In what follows, we start by showing that the Kullback--Leibler divergence, while relatively practical from a computational perspective, 
can be small even when the approximate point estimates and uncertainties are far (sometimes arbitrarily far) from the exact values.
By contrast, we show that closeness in $p$-Wasserstein distance (with $p=1$ or $2$) implies closeness in relevant point estimates and uncertainties. 
Unfortunately, though, the Wasserstein distance is challenging to use in practice. 
To address this shortcoming, we introduce the $(p,\nu)$-Fisher distance, which generalizes a number of existing distances
in a range of literatures~\citep{Johnson:2004a,Johnson:2004,Hyvarinen:2005,Bolley:2012,Sriperumbudur:2017,Huggins:2017a,Huggins:2017b,Campbell:2017,Campbell:2018}.
We extend and synthesize the results of  \citet{Bolley:2012}  and  \citet{Huggins:2017a}
to show that the $(p,\nu)$-Fisher distance provides an upper bound on the $p$-Wasserstein distance in many cases of interest. 
We illustrate that the $(p,\nu)$-Fisher distance avoids many of the pitfalls of the Kullback--Leibler divergence. 
We also show that the $(p,\nu)$-Fisher is more practical to calculate than the Wasserstein distance.

We demonstrate the practicality of our proposed $(p,\nu)$-Fisher distance by using it to analyze 
two scalable Bayesian approximation methods: the Laplace approximation~\citep{Schervish:1995} and Bayesian coresets~\citep{Campbell:2017,Campbell:2018}.
First, we derive computable bounds on the $p$-Wasserstein distance between the exact posterior and the Laplace approximation 
for Bayesian models with log posterior densities that are strongly convex and have bounded third derivatives. 
As a corollary we provide a bound on the convergence rate of the Laplace approximation to the exact posterior in $p$-Wasserstein distance. 
Second, we consider the accuracy of using a \emph{coreset}, which is a small, weighted subset of data, to approximate the likelihood. 
\citet{Campbell:2017} created Bayesian coresets designed to provide high-quality posterior approximations. 
A coreset approximation can be computed more quickly than the full likelihood across all data points and thus in turn can be used as a 
pre-processing step to speed up standard inference methods such as Markov chain Monte Carlo.
Our results together with those from \citet{Campbell:2018} 
imply that the $p$-Wasserstein distance between the exact and coreset posteriors decreases exponentially in the size of the coreset. 

\section{Preliminaries} 
\label{sec:preliminaries}

Let $\data = (z_{1}, \dots, z_{n})$ denote our observed data, and let $\param \in \reals^{d}$ denote our parameter
vector of interest. A Bayesian model consists of a prior measure $\priordist(\dee\param)$ and a likelihood $\lik(\data; \param)$. Together the prior and likelihood define a joint distribution over the data and parameters. The Bayesian posterior distribution is the conditional in $\param$. To write this conditional, we define the log likelihood $\loglik(\param) = \log \lik(\data; \param)$ and marginal likelihood, or evidence, $\marginallik = \int \exp\{\loglik(\param)\}\priordist(\param)\dee\param$. Then the posterior is
\[
\postdist(\dee\param) = \frac{e^{\loglik(\param)}\priordist(\dee\param)}{\marginallik}. \label{eq:posterior}
\]
Since in Bayesian analysis the data are fixed and we will treat the data as constant for the remainder of the paper, 
we have suppressed the dependence on $Z$ in our notation. %

Typically practitioners report summaries of the posterior distribution in the form of point estimates and uncertainties; we introduce a number of relevant summaries here. 
For some distribution $\eta$ on $\reals^{d}$, let $\Param \dist \eta$.
Let $\mean{\eta} = \EE(\Param)$ denote the mean of $\Param$. 
Let $\Sigma_{\eta} = \EE\{(\Param - \meanvec{\eta})(\Param - \meanvec{\eta})^{\top}\}$ denote the covariance of $\Param$, and let 
$\std{\eta,i} = \Sigma_{\eta,ii}^{1/2}$~$(i=1,\ldots,d)$ denote the standard deviation of the $i$th component of $\Param$.
An alternative measure of uncertainty is the mean absolute deviation, $\MAD{\eta,i} = \EE(|\Param_{i} - \mean{\eta,i}|)$~$(i=1,\ldots,d)$. 
To construct medians and quantiles, we define $I_{\eta,i,a,b} = \EE\{\ind_{[a,b]}(\Param_{i})\}$ ($i=1,\dots,d; \infty \le a < b \le \infty$), 
where $\ind_{A}(\param)$ equals 1 when $\param \in A$ and $0$ otherwise. When $d = 1$, we drop the index $i$ from the subscript. 
We use standard asymptotic notation: $f = O(g)$ if and only if $\limsup f/g < \infty$ and $f = \Theta(g)$ if and only if $\limsup f/g < \infty$ and $\liminf f/g > 0$. 
Although we also use $\Theta$ to denote a random variable, the meaning will be clear from context. 

\section{Kullback--Leibler divergence}
\label{sec:KL}

Let $\estdist$ be any approximation to the posterior $\postdist$; we take $\estdist$ to be a Borel probability measure. 
In order to bound the error in approximate summaries derived from $\estdist$, we consider an intermediate step of first bounding some notion of divergence between $\estdist$ and $\postdist$. 
A choice notion of divergence to consider is the Kullback--Leibler divergence since one of the most widely-used posterior 
approximation methods, \emph{variational inference}, works by minimizing the Kullback--Leibler divergence (in a particular direction) 
over a tractable family $\varfamily$ of potential approximation distributions~\citep{Blei:2017}:
\[ \label{eq:klmin}
\estdist = \argmin_{\vardist \in \varfamily} \kl{\vardist}{\postdist},
\]
where
\[
\kl{\vardist}{\postdist} = \int \log \der{\vardist}{\postdist}(\param) \vardist(\dee\param).
\]
Part of what makes $\varfamily$ tractable is typically that we can compute the expectations needed to solve \cref{eq:klmin}. 
In particular, for any $\vardist \in \varfamily$, we are usually able to efficiently calculate relevant summaries, such as those in \cref{sec:preliminaries}, either analytically
or using independent and identically distributed samples from $\vardist$.  
By contrast, we cannot expect easy access to the moments of $\postdist$. 
Therefore, we might ask if a small Kullback--Leibler divergence implies that these approximate summaries have small error.

To emphasize the generality of our results beyond Bayesian inference, 
in what follows we let $\mainmeas$ and $\estmeas$ denote two Borel probability measures.
We will typically take $\mainmeas = \postdist$ and $\estmeas = \estdist$ in the Bayesian case, 
so we imagine that we have access to the summaries of $\estmeas$ but not $\mainmeas$. 
First, we notice that if the Kullback--Leibler divergence between $\estmeas$ and $\mainmeas$ is small, then $\estmeas$ can 
provide good mean and quantile estimates.
\begin{proposition} \label{prop:KL-divergence-working}
If $\delta = \kl{\estmeas}{\mainmeas} < 1$, then $(\mean{\estmeas} - \mean{\mainmeas})^{2} \le (\var{\estmeas} + \var{\mainmeas})\delta/(1 - \delta)$
and $|I_{\estmeas,a,b} - I_{\mainmeas,a,b}| \le (\delta/2)^{1/2}$. 
\end{proposition}

In order for small Kullback--Leibler divergence to imply good mean estimates, \cref{prop:KL-divergence-working} 
also requires the variances $\var{\estmeas}$ and $\var{\mainmeas}$ not to be too large; we examine this assumption  below. 
\cref{prop:KL-divergence-working} guarantees good credible intervals in the sense that any $\hat\alpha$-confidence interval for $\estmeas$ 
will be an $\alpha$-confidence interval for $\mainmeas$, where $\alpha \in (\hat\alpha - (\delta/2)^{1/2}, \hat\alpha + (\delta/2)^{1/2})$. 
On the other hand, as we see in the next result, even when the Kullback--Leibler divergence between $\estmeas$ and $\mainmeas$ is small the 
variance estimate provided by $\estmeas$ can be arbitrarily bad. Moreover, if the Kullback--Leibler divergence is moderately sized, the mean estimates may be very far off.
\begin{proposition} \label{prop:KL-divergence-problems}
(A) For any $\delta > 0$ there exist Gaussian distributions $\estmeas$ and $\mainmeas$ such that $\kl{\estmeas}{\mainmeas} = \delta$,
$(\mean{\mainmeas} - \mean{\estmeas})^{2} = \var{\estmeas}\{\exp(2\delta) - 1\}$, and $\var{\estmeas} = \exp(-2\delta)\var{\mainmeas}$.

(B) For any $t \in (1, \infty]$, there exist mean-zero unimodal distributions $\estmeas$ and $\mainmeas$ such that 
$\kl{\estmeas}{\mainmeas} < 0{\cdot}802$ but $\var{\mainmeas} \ge t\var{\estmeas}$. 

\end{proposition}
\begin{remark}
In part (B), the distributions used are very simple: $\estmeas$ is a standard Gaussian and $\mainmeas$ is 
a standard $t$-distribution with $h \ge 2$ degrees of freedom. 
Numerical computations suggest that, the constant $0{\cdot}802$ can be replaced by $0{\cdot}12$.
\end{remark}

Part (A) of \cref{prop:KL-divergence-problems} shows that, for example, if $\kl{\estmeas}{\mainmeas} = 5$ then
the mean estimate may be off by more than $148\std{\estmeas}$. 
Since $\std{\estmeas}$ provides a natural unit of uncertainty about the parameter value, 
we see that a moderate Kullback--Leibler divergence can correspond to a very large error in the mean estimate. 
Note in particular that since variational methods typically optimize over a constrained set of tractable distributions such 
as products of exponential families~\citep{Ranganath:2014,Kucukelbir:2015,Blei:2017},
moderate Kullback--Leibler values are expected in many applications.
Part (B) shows that unless the Kullback--Leibler divergence is very small, the posterior may have arbitrarily large 
variance no matter the approximate variance observed from $\estmeas$.
Therefore, the Kullback--Leibler divergence is able to capture the quality of mean and uncertainty 
estimates only in limited circumstances and cannot be relied upon to 
capture the quality of variance estimates.  %

\section{Wasserstein distance}
\label{sec:Wasserstein}

A suggestion for an alternative divergence to consider is provided by the theory supporting Markov chain Monte Carlo, where the \emph{Wasserstein distance} is widely used~\citep{Joulin:2010,Madras:2010,Hairer:2014,Rudolf:2015,Durmus:2016b,Vollmer:2016,Durmus:2017,Cheng:2017,Mangoubi:2017,Cheng:2018,Fang:2018,BouRabee:2018}.
Wasserstein distance has also been adopted for asymptotic analysis in the large data limit~\citep{Minsker:2017}.
Let $\couplings{\mainmeas}{\estmeas}$ denote the set of couplings between $\mainmeas$ and $\estmeas$.
That is, $\couplings{\mainmeas}{\estmeas}$ is the set of Borel measures $\coupling$ on $\reals^{d} \times \reals^{d}$ such 
that $\coupling$ has marginal distributions $\mainmeas$ and $\estmeas$:
$\mainmeas = \coupling(\cdot, \reals^{d})$ and $\estmeas = \coupling(\reals^{d}, \cdot)$. 
The $p$-Wasserstein distance between $\mainmeas$ and $\estmeas$ is given by~\citep[Def.~6.1]{Villani:2009}
\[
\pwassSimple{p}{\mainmeas}{\estmeas} 
= \inf_{\coupling \in \couplings{\mainmeas}{\estmeas}} \left\{ \int \staticnorm{\param - \hat{\param}}_{2}^{p} \coupling(\dee \param, \dee \hat{\param}) \right\}^{1/p}.
\]
By Jensen's inequality, 
\[
\pwassSimple{p'}{\mainmeas}{\estmeas} 
	&\le \pwassSimple{p}{\mainmeas}{\estmeas} 
	\qquad (1 \le p' \le p < \infty). \label{eq:p-Wasserstein-ordering}
\]
The well-known dual form of the 1-Wasserstein distance can be convenient to work with and is particularly interpretable.
For a function $\phi : \reals^{d} \to \reals$, let 
$\norm{\phi}_{L} = \sup_{\param \ne \param'} |\phi(\param) - \phi(\param')|/\twonorm{\param - \param'}$ 
denote its Lipschitz norm. 
Let $\mainmeas(\phi) = \int \phi(\param)\mainmeas(\dee\param)$ denote the expectation of the integrable function $\phi$
with respect to the measure $\mainmeas$.
Then~\citep[Rmk.~6.5]{Villani:2009}
\[
\pwassSimple{1}{\mainmeas}{\estmeas}  = \sup_{\phi \st \norm{\phi}_{L} \le 1} |\mainmeas(\phi) - \estmeas(\phi)|.  \label{eq:1-Wasserstein-dual}
\]
Eqs.~\eqref{eq:p-Wasserstein-ordering} and \eqref{eq:1-Wasserstein-dual} together imply that for any $p \ge 1$, 
if $\pwassSimple{p}{\mainmeas}{\estmeas} \le \veps$, then for any $L$-Lipschitz function $\phi$, 
$|\mainmeas(\phi) - \estmeas(\phi)| \le L\veps$. It follows that we can obtain bounds on the error of our approximate summaries.
\begin{theorem} \label{thm:Wasserstein-moment-bounds}
(A) If $\pwassSimple{1}{\mainmeas}{\estmeas} \le \veps$ then 
$
\twonorm{\mean{\mainmeas} - \mean{\estmeas}} \le \veps 
$
and
$
|\MAD{\mainmeas,i} - \MAD{\estmeas,i}| \le 2\veps %
$ $(i=1,\dots,d)$. 
If in addition $\mainmeas$ has Lebesgue density bounded by $c < \infty$, then
\[
|I_{\mainmeas,i,a,b} - I_{\estmeas,i,a,b}| 
	&\le 2(2 c\,\veps)^{1/2} 
	 \qquad (i=1,\dots,d; \infty \le a < b \le \infty). \label{eq:credible-interval-error-bound}
\]
(B) If $\pwassSimple{2}{\mainmeas}{\estmeas} \le \veps$, then 
$\pwassSimple{1}{\mainmeas}{\estmeas} \le \veps$,
\[
|\std{\mainmeas,i} - \std{\estmeas,i}| \le \frac{1}{2} \left(2^{1/2}+6^{1/2}\right)\veps \le 2\veps \qquad (i=1,\dots,d), %
\]
and
\[
\twonorm{\Sigma_{\mainmeas} - \Sigma_{\estmeas}}
&\le 2^{3/2}\min(\twonorm{\Sigma_{\mainmeas}}^{1/2},\twonorm{\Sigma_{\estmeas}}^{1/2})\veps + (1 + 3 \times 2^{1/2})\veps^{2} \\
&< 3\min(\twonorm{\Sigma_{\mainmeas}}^{1/2},\twonorm{\Sigma_{\estmeas}}^{1/2})\veps + 5{\cdot}25\veps^{2}.
\]
\end{theorem}

By taking $\mainmeas = \postdist$ and $\estmeas = \estdist$,
\cref{thm:Wasserstein-moment-bounds} shows that
the Wasserstein distance can be used to bound the error of estimates of the posterior mean,
covariance matrix, standard deviation, mean absolute deviation, and credible intervals. 
A weakness of \cref{thm:Wasserstein-moment-bounds}, which we will remedy shortly, 
is that the bound on $|I_{\mainmeas,i,a,b} - I_{\estmeas,i,a,b}|$ requires $\mainmeas$ 
to not be too peaked, as otherwise the constant $c$ will be large. 
However, in the large data (i.e., large $n$) limit when a Bernstein--Von Mises theorem applies, 
we expect $c = \Theta(n^{1/2})$.
Hence $\veps = O(n^{-1/2})$ for \cref{eq:credible-interval-error-bound} to be nontrivial.

\begin{remark}
If $\mainmeas$ and $\estmeas$ are univariate Gaussian distributions then
\[
\pwassSimple{2}{\mainmeas}{\estmeas}^{2} 
&= (\mean{\mainmeas} - \mean{\estmeas})^{2} + (\std{\mainmeas} - \std{\estmeas})^{2}.
\]
So if $\pwassSimple{2}{\mainmeas}{\estmeas} \le \veps$, we can conclude
that $|\mean{\mainmeas} - \mean{\estmeas}| \le \veps$ and $|\std{\mainmeas} - \std{\estmeas}| \le \veps$.
Thus, the error bounds given for 2-Wasserstein are tight for the mean and tight up 
to a factor of $(2^{1/2}+6^{1/2})/2$ for the standard deviation. 
If $\std{\mainmeas} = \std{\estmeas}$, then $\pwassSimple{1}{\mainmeas}{\estmeas} = |\mean{\mainmeas} - \mean{\estmeas}|$,
so the mean moment error bound given for 1-Wasserstein is tight. 
\end{remark}

We have seen that Wasserstein distance provides exactly the guarantees on summary approximation error that we were looking for, but it poses a number of computational challenges. For one, the intractable normalizing constant (that is, the marginal likelihood $\marginallik$) is still present in the exact Bayesian posterior $\mainmeas = \postdist$; by contrast, the Kullback--Leibler optimization problem in \cref{eq:klmin} can be solved by using the unnormalized version of the posterior instead~\citep[Ch.~10]{Bishop:2006}. Moreover, the $\inf$ (or $\sup$ in the dual formulation) poses an added challenge. For these reasons, Wasserstein distance, unlike the Kullback--Leibler divergence, is rarely used as an optimization objective, with \citet{Srivastava:2015} a notable exception that takes advantage of the measures under consideration being discrete.

\section{Wasserstein distance bounds via the $(p, \nu)$-Fisher norm}
\label{sec:Fisher}

We introduce a new statistical distance, which we call the \emph{$(p,\nu)$-Fisher distance},
as an alternative that is more computationally tractable. We show below that the
$(p,\nu)$-Fisher distance implies a bound on the $p$-Wasserstein distance in many cases of interest. 
And therefore the $(p,\nu)$-Fisher distance in turn implies a bound on the error of approximate posterior summaries by \cref{thm:Wasserstein-moment-bounds}. %
For a Borel measure $\nu$, let $\Lparg{p}{\nu}$ denote the space of functions that are $p$-integrable with respect to $\nu$: 
$\phi \in \Lparg{p}{\nu}  \iff  \Lpnormarg{p}{\nu}{\phi} = (\int \phi(\param)^{p}\nu(\dee \param))^{1/p} < \infty$. 
Let $\mainPE = -\log\textder{\mainmeas}{\param}$ and $\estPE = -\log\textder{\estmeas}{\param}$ denote the 
potential energy functions associated with, respectively,  $\mainmeas$ and $\estmeas$.
\begin{definition}
The \emph{$(p,\nu)$-Fisher distance} is given by
\[
\fd{p,\nu}{\mainmeas}{\estmeas} 
&= \Lpnormarg{p}{\nu}{\statictwonorm{\grad\mainPE - \grad\estPE}} 
= \left\{\int \statictwonorm{\grad\mainPE(\param) - \grad\estPE(\param)}^{p} \nu(\dee\param)\right\}^{1/p}.
 \label{eq:fisher-distance} 
\]
\end{definition}
The special case $I(\mainmeas|\estmeas) = \fd{p,\mainmeas}{\mainmeas}{\estmeas}$ is known by many names, 
including the Fisher divergence \citep{Sriperumbudur:2017} and the Fisher information of $\mainmeas$ 
with respect to $\estmeas$~\citep{Bolley:2012}.
The Fisher divergence has been used to prove central limit theorems \citep{Johnson:2004a,Johnson:2004} 
and as an objective for density estimation \citep{Sriperumbudur:2017,Hyvarinen:2005}.
Special cases of the $(p,\nu)$-Fisher distance have also been used in a Bayesian context both for analyzing approximation quality \citep{Huggins:2017a,Huggins:2017b} and as an objective function for approximate inference~\citep{Campbell:2017,Campbell:2018}.
We will discuss some of these applications in detail in \cref{sec:applications}.

In the Bayesian posterior case where $\mainmeas = \postdist$, we note that the computationally intractable posterior normalizer $\marginallik$ 
is constant in $\param$ and therefore vanishes in the gradient $\grad\mainPE$. 
Hence the $(p,\nu)$-Fisher distance avoids the principal computational challenges of the Wasserstein distance. 
Our next results show that $(p,\nu)$-Fisher distance also bounds Wasserstein distance.
Specifically, for well-behaved densities, the $p$-Wasserstein distance between $\mainmeas$ and $\estmeas$
is bounded by a multiple of the $(p,\mainmeas)$-Fisher distance.
\begin{theorem} \label{thm:p-Wasserstein-sde-error}
Assume $\mainPE$ and $\estPE$ are twice continuously differentiable
and that for some $\alpha > 0$, $\estPE$ is $\alpha$-strongly convex:
\[
\estPE(\param') \ge \estPE(\param)+\grad \estPE(\param)^{\top}(\param' - \param) + (\alpha/2)\statictwonorm{\param - \param'}^{2}
\quad (\param,\param' \in \reals^{d}).
\]
Then for $p \in \{1,2\}$, 
\[
\pwassSimple{p}{\mainmeas}{\estmeas} 
\le \alpha^{-1}\fd{p,\mainmeas}{\mainmeas}{\estmeas}.  \label{eq:p-Wasserstein-sde-error}
\]
\end{theorem}
\bprf
We prove the $p=1$ case in the Appendix. 
The $p=2$ case follows from \citet[Lemma 3.3 and Proposition 3.10]{Bolley:2012}.
\eprf
Requiring $\estPE$ to be $\alpha$-strongly convex is a widespread assumption in analyses of 
Markov chain Monte Carlo algorithms~\citep{Durmus:2016b,Vollmer:2016,Cheng:2017,Mangoubi:2017}.
However, it is a strong assumption which can be weakened to, essentially, only assuming strong
convexity of $\estPE$ outside some compact set.
\begin{theorem} \label{thm:p-Wasserstein-sde-error-general}
Fix $p \in \{1,2\}$. 
Assume $\mainPE$ and $\estPE$ are twice continuously differentiable
and for some constants $K > 0$ and $R \ge 0$, $\hessian \estPE(\param) \succeq K I_{d}$ for all $\twonorm{\param} \ge R$.
If $p=1$, further assume that each continuous function $\phi$ is $\estmeas$-integrable if it is $\mainmeas$-integrable. 
Then for $\alpha = \alpha(p, \estPE, R, K)$ but independent of $\mainPE$,
\[
\pwassSimple{p}{\mainmeas}{\estmeas} 
\le \alpha^{-1}\fd{p,\mainmeas}{\mainmeas}{\estmeas}.  \label{eq:p-Wasserstein-sde-error-general}
\]
\end{theorem}
\bprf
The $p=1$ result follows from \citet[Thm.~3.4]{Huggins:2017a} and \citet[Cor.~2]{Eberle:2015}.
The $p=2$ result follows from \citet[Proposition 3.4 and 3.10]{Bolley:2012}.
\eprf
\begin{remark}
In the $p = 1$ case, the condition on $\estPE$ can be further weakened to a condition \citet{Gorham:2016b} 
call \emph{distant dissipativity}~\citep[see also][]{Eberle:2015,Huggins:2017a}.
\end{remark}
\begin{remark}
Although $\alpha$ depends on $\estPE$, it does so only through a limited number of properties.
For example, in the $p=2$ case, $\alpha$ depends on the minima and maxima of $\estPE$
on the ball of center 0 and radius $(1+\veps)R$ for any choice of $\veps > 0$.
\end{remark}

A limitation of Theorems \ref{thm:p-Wasserstein-sde-error} and \ref{thm:p-Wasserstein-sde-error-general}
is that they bound the $p$-Wasserstein distance
in terms of the $(p,\nu)$-Fisher distance only when $\nu = \mainmeas$. 
However, we would like the flexibility to handle more general choices of $\nu$. In particular, an integral
with respect to $\mainmeas = \postdist$ is typically computationally intractable in the Bayesian case, so
we wish to consider more tractable choices for $\nu$.

Take any Borel probability measures $\xi$ and $\nu$ with $\xi$ absolutely continuous with respect to $\nu$ (i.e., $\xi \abscont \nu$);
in this case 
 the \emph{$\chi^{2}$-divergence}~\citep{Csiszar:1967} is defined as
$
{\dchi{\xi}{\nu} = \int \textder{\xi}{\nu}(\param) \xi(\dee\param) - 1}.
$
\begin{corollary} \label{cor:Wasserstein-error}
Assume the hypotheses of \cref{thm:p-Wasserstein-sde-error} hold. 
Let $B_{1}(\mainmeas, \nu) = (1 + \dchi{\mainmeas}{\nu})^{1/2}$ 
and $B_{2}(\mainmeas, \nu) =  \infnorm{\textder{\mainmeas}{\nu}}^{1/2}$. 
Then for any probability measure $\nu$ such that $\mainmeas \abscont \nu$, 
\[
\pwassSimple{p}{\mainmeas}{\estmeas} 
\le  \alpha^{-1}B_{p}(\mainmeas, \nu)\fd{2,\nu}{\mainmeas}{\estmeas}  \quad (p=1,2) \label{eq:p-Wasserstein-error}
\]
\begin{remark}
For bounding both 1-Wasserstein and 2-Wasserstein distance, \cref{cor:Wasserstein-error} relies on the $(2,\nu)$-Fisher distance. 
\end{remark}
\end{corollary}

\subsection{Tightness of the bounds}

The Wasserstein bounds provided by \cref{cor:Wasserstein-error} eliminate the possibility of dangerous situations as
in \cref{prop:KL-divergence-problems}(B), where the Kullback--Leibler divergence between $\estmeas$ and $\mainmeas$ 
was finite but the distribution of interest $\mainmeas$ had arbitrarily large or infinite variance. 
But it remains to show that the bounds are tight enough for practical use. 
To investigate this question we consider the two settings from \cref{prop:KL-divergence-problems}
and for simplicity focus on the $p=2$ case.
We start by considering the Gaussian setting, as in \cref{prop:KL-divergence-problems}(A).

\begin{figure}
\centering
\includegraphics[width=0.6\textwidth]{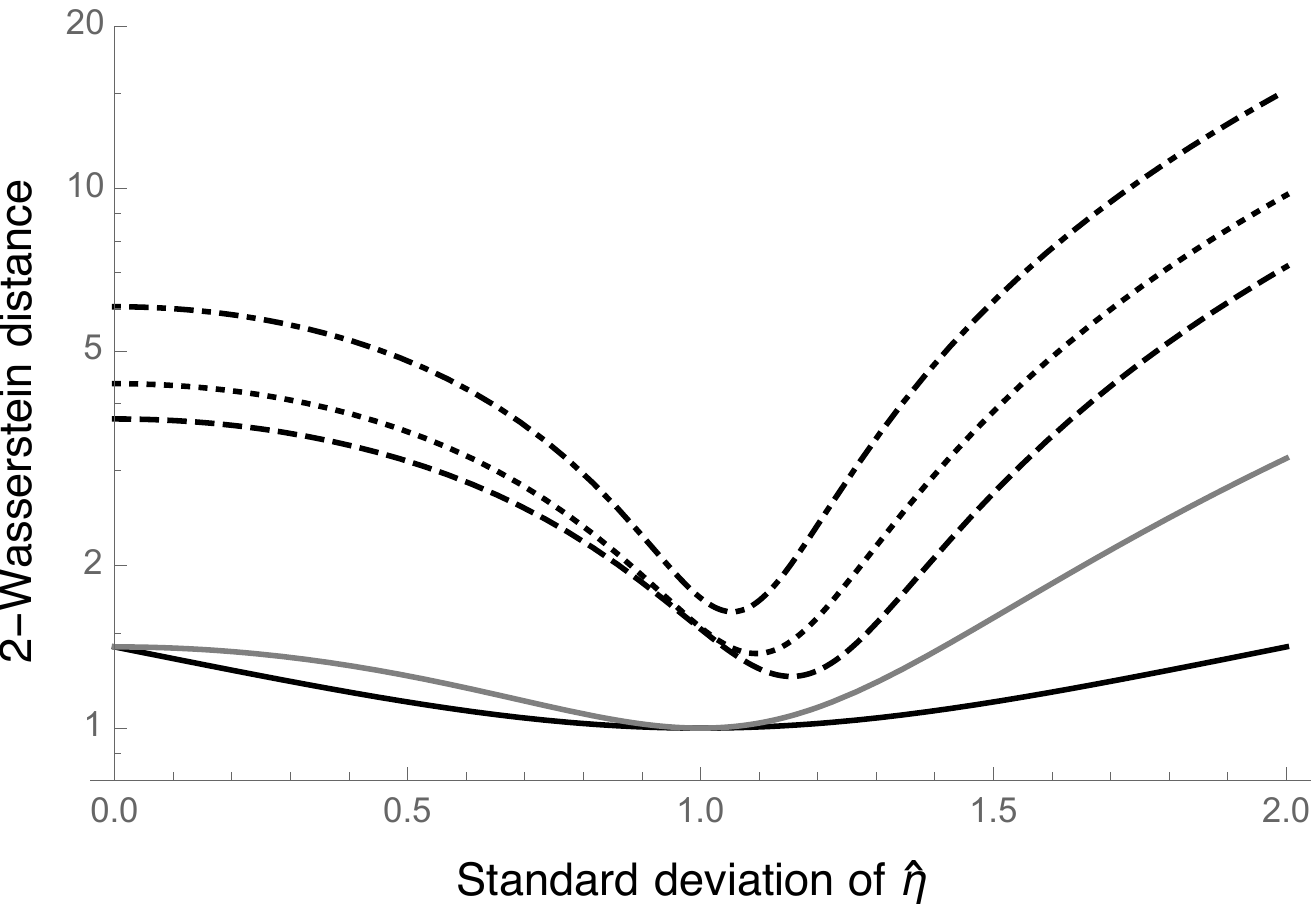}
\caption{When $\std{\mainmeas} = 1$ and $\Delta\mean{} = -1$, the value of $\pwassSimple{2}{\mainmeas}{\estmeas}$ (solid black),
the bound from \cref{eq:nu-equal-post-bound} (solid gray), and bounds from \cref{eq:nu-normal-bound} when $\eps = 1$ and $\rho = 2$ (dashes), 
$\rho = 4$ (dots), and $\rho = 8$ (dot-dashes). }
\label{fig:2-Wasserstein-bounds}
\end{figure}

\begin{proposition} \label{prop:normal-bound-tightness}
Let $\Delta\mean{} = \mean{\mainmeas} - \mean{\estmeas}$,  $\Delta\std{} = \std{\mainmeas} - \std{\estmeas}$, 
and $r = \std{\estmeas}/\std{\mainmeas}$. 
If $\mainmeas$ and $\estmeas$ are Gaussian, then $\pwassSimple{2}{\mainmeas}{\estmeas}^{2} = (\Delta\mean{})^{2} + (\Delta\std{})^{2}$ while \cref{cor:Wasserstein-error} implies that when $\nu = \distNorm(\mean{\mainmeas} + \eps, \rho \var{\mainmeas})$ for some  $\eps \in \reals$ and $\rho > 1$, 
\[
\pwassSimple{2}{\mainmeas}{\estmeas}^{2} 
\le  C \left\{ r^{2}(r^{2} + 1)\eps^{2} + (1- r^{2})(\eps - \Delta\mean{})^{2} + r^{2}(\Delta\mean{})^{2} + \rho (r + 1)^{2} (\Delta\std{})^{2} \right\}, \label{eq:nu-normal-bound}
\]
where $C =\rho^{1/2}\exp[\epsilon^2 / \{2 (\rho-1)\var{\mainmeas}\}]$.
In particular, when $\nu = \mainmeas$,  \cref{cor:Wasserstein-error} implies that 
\[
\pwassSimple{2}{\mainmeas}{\estmeas}^{2}  \le \Delta^{2} +  (1 + r)^{2}(\Delta\std{})^{2}. \label{eq:nu-equal-post-bound}
\]
\end{proposition}
\cref{eq:nu-equal-post-bound} shows that in the ideal case of $\nu = \mainmeas$, the 2-Wasserstein bound is quite tight. However, 
the bound in \cref{eq:nu-normal-bound} is more difficult to interpret. 
\cref{fig:2-Wasserstein-bounds} provides some additional insight by considering the behavior of the bounds
for fixed $\std{\mainmeas}, \Delta\mean{}$, and $\eps$. 
The figure confirms that the $\nu = \mainmeas$ bound is reasonably accurate while the bounds when $\nu \ne \mainmeas$ are looser.
When $\rho$ is small the bound is tighter but its minimum is farther from the true optimum. 
When $\rho$ is larger the bound is looser but the optimum approaches the correct value of one. 
Thus, when $\mean{\nu}$ is incorrect (that is, $\eps \ne 0$), there appears to be a tightness-bias tradeoff when selecting $\var{\nu}$ (that is, $\rho$). 

\begin{figure}
\centering
\includegraphics[width=0.6\textwidth]{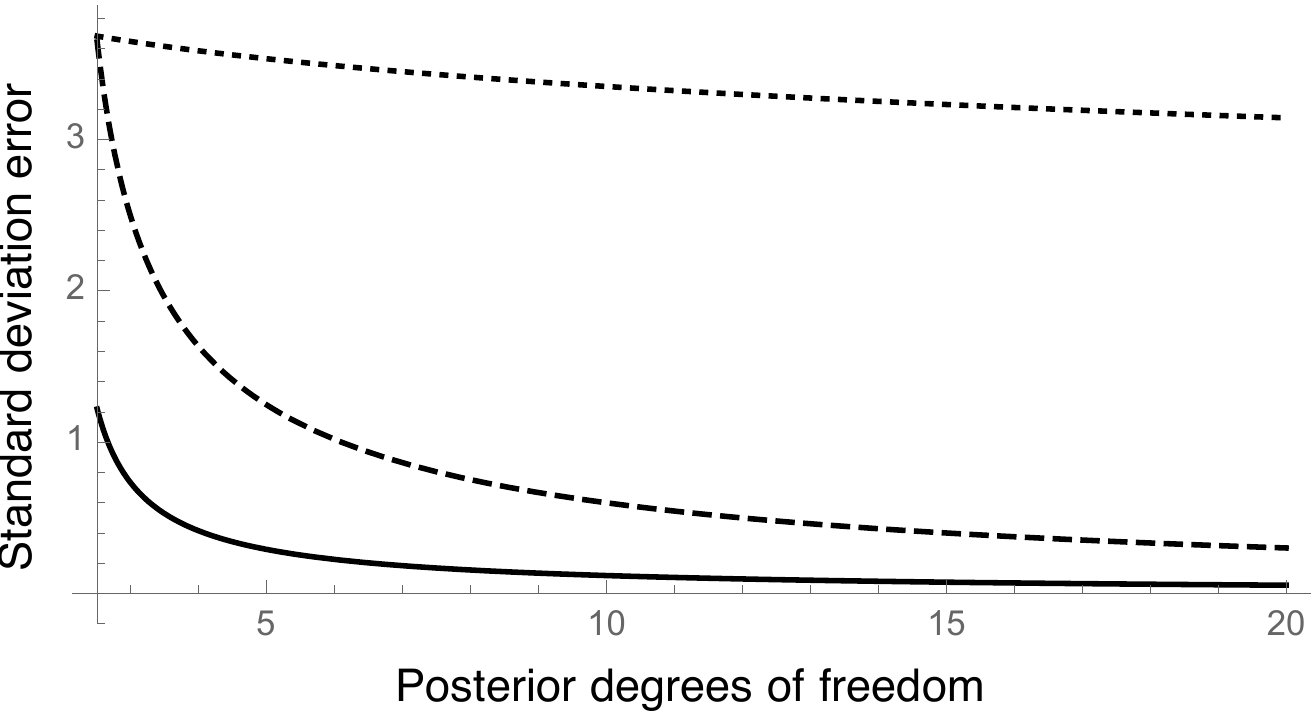}
\caption{The true standard deviation error $|\std{\mainmeas} - \std{\estmeas}|$ as a function of $h$ (solid) along with the 2-Wasserstein bounds on the error
when using $\nu = \mainmeas$ (dashes) and $\nu = \mainmeas_{h_{0}}$ when $h_{0} = 2.5$ (dots).}
\label{fig:std-dev-error}
\end{figure}

In \cref{prop:KL-divergence-problems}(B), we took $\estmeas = \distNorm(0,1)$ and $\mainmeas = \mcT_{h}$, a standard $t$-distribution 
with $h \ge 2$ degrees of freedom.
The distribution $\mainmeas$ could arise as a posterior, for example, by placing an normal-inverse gamma prior on the mean and variance of a normal observation model 
and integrating out the variance. 
Equivalently, we could place a $t$-distribution prior with $h_{0} = h - n$ degrees of freedom on the mean of a normal observation model.
\begin{proposition} \label{prop:t-dist-bound-tightness}
If $\estmeas = \distNorm(0,1)$ and $\mainmeas = \mcT_{h}$, then 
\cref{cor:Wasserstein-error} implies that when $\nu =  \mainmeas$,
\[
\pwassSimple{2}{\mainmeas}{\estmeas}^{2} \le 10/(h^{2} + h - 6). \label{eq:t-dist-wasserstein}
\]
\end{proposition}
Figure~\ref{fig:std-dev-error} shows that bound on the standard deviation implied by \cref{eq:t-dist-wasserstein} is 
reasonably tight, while the bound obtained when $\nu = \mcT_{h_{0}}$ (using a computer algebra system) is somewhat looser.

\subsection{Strong convergence and $(p,\nu)$-Fisher distance}

While convergence in Wasserstein distance implies weak convergence, it does not imply strong convergence.
The $(2,\nu)$-Fisher distance, on the other hand, does imply convergence in total variation distance 
$d_{\mathrm{TV}}(\mainmeas, \estmeas) = \sup_{A \subseteq \reals^{d}} |\mainmeas(A) - \estmeas(A)|$
and hence strong convergence. 

\begin{proposition} \label{prop:Fisher-TV-bound}
If \cref{eq:p-Wasserstein-sde-error-general} holds then 
\[
d_{\mathrm{TV}}(\mainmeas, \estmeas)  \le (2\alpha)^{-1/2}  \infnorm{\textder{\mainmeas}{\nu}}^{1/2}\fd{2,\nu}{\mainmeas}{\estmeas}. \label{eq:Fisher-TV-bound}
\]
In particular, $|I_{\mainmeas,i,a,b} - I_{\estmeas,i,a,b}|~(i=1,\dots,d; \infty \le a < b \le \infty)$ is 
upper bounded by the left-hand side of \cref{eq:Fisher-TV-bound}.
\end{proposition}

\cref{prop:Fisher-TV-bound} remedies a shortcoming of \cref{thm:Wasserstein-moment-bounds} discussed earlier: 
that the bound on $|I_{\mainmeas,i,a,b} - I_{\estmeas,i,a,b}|$ requires a good bound on the density of $\mainmeas$. 
\citet{Johnson:2004a}, \citet{Johnson:2004}, and \citet{Ley:2013} provide similar bounds in the one-dimensional case
for  certain integral probability measures such as the total variation and Kolmogorov distances. 

\section{Applications}
\label{sec:applications}

We consider two applications of the $(p,\nu)$-Fisher distance for controlling the Wasserstein error 
of approximate Bayesian inference methods. 

\subsection{Laplace approximation}

For a function $\phi : \reals^{d} \to \reals$, define its matrix of second partial derivatives 
$(\derivop{2} \phi)(\param)_{ij} = \partial_{i}\partial_{j}\phi(\param)~{(i=1,\dots,d}; {j=1,\dots,d)}$.
Also, with $\postdensity$ denoting the density of $\postdist$ with respect to Lebesgue measure on $\reals^{d}$, 
let 
$
\map = \argmax_{\param} \log\{\postdensity(\param)\} %
$
denote the maximum a posteriori parameter estimate. 
Denote the Hessian matrix at $\map$ by ${\maphessian = (\hessian\log\postdensity)(\map)}$. 
Then the Laplace approximation for $\postdist$ is the Gaussian distribution
\[
\approxdist_{\mathrm{Laplace}}(\dee \param) = (2\pi)^{-d/2}|\maphessian|^{1/2}\exp\{0{\cdot}5(\param - \map)^{\top}\maphessian(\param - \map)\}\dee\param. 
\]

\begin{proposition}[Non-asymptotic Laplace approximation error] \label{prop:nonasymptotic-Laplace}
Assume that $-\log\postdensity$ is three times continuously differentiable and $\alpha$-strongly convex, and that
\[
\sup_{\param}\sum_{i=1}^{d}\twonorm{\{\hessian\partial_{j}\log\postdensity\}(\param)}^{2} \le \mathsf{M}^{2}.
\]
Let $\lambda$ denote the eigenvalues of $(\maphessian)^{-1}$, $\mathsf{L}_{1}(\lambda) = \onenorm{\lambda}$, 
and $\mathsf{L}_{2}(\lambda) = (\onenorm{\lambda}^{2} + 2\twonorm{\lambda}^{2})^{1/2}$. 
Then 
\[
\pwassSimple{p}{\approxdist_{\mathrm{Laplace}}}{\postdist}  \le \alpha^{-1}\mathsf{L}_{p}(\lambda) \mathsf{M} \quad (p = 1,2).
\]
\end{proposition}

In many settings the bounds provided by \cref{prop:nonasymptotic-Laplace} are computable.
The key is to show that $-\log\postdensity$ is strongly convex and that the third derivatives of $\log\postdensity$ are uniformly bounded. 
For typical models with $n$ observations that have Berstein--Von Mises behavior, $\mathsf{M} = \Theta(n)$ while $\mathsf{L}_{p} = \Theta(n^{-1})$. 
Hence, the bound on $\pwassSimple{p}{\approxdist_{\mathrm{Laplace}}}{\postdist}$
is determined by $\alpha$, the strong convexity constant of $-\log\postdensity$.
Unfortunately, for many models of interest $\alpha = O(1)$ because the negative log-likelihood is convex but not strongly convex.
To see that this is not likely to be an issue in practice we can consider an asymptotic analysis under assumptions 
that are standard in the Laplace approximation literature~\citep{Tierney:1989,Kass:1990,Schervish:1995,Small:2010}. 
To state our result, we consider a sequence of absolutely continuous posterior distributions $(\postdist_{n})_{n=1}^{\infty}$
and define the normalized log posterior densities $\logpostn = n^{-1}\log \textder{\postdist_{n}}{\param}~(n = 1,\dots)$, 
where $n$ can be thought of as the number of observations available. 
Let $\mapn = \argmax_{\param}\logpostn(\param)$ and  $\maphessiann = (\hessian\logpostn)(\mapn)$. 

\begin{proposition}[Asymptotic Laplace approximation error] \label{prop:asymptotic-Laplace}
Assume that $\logpostn$ is three times continuously differentiable, 
$\norm{\maphessiann} \ge \alpha^{-1}$, 
$\sup_{\param}\sum_{i=1}^{d}\twonorm{\{\hessian\partial_{j}\logpostn\}(\param)}^{2} \le \mathsf{M}^{2}$,
and 
\[
\left\{\int\twonorm{\param - \mapn}^{2p} \postdist_{n}(\dee \param)\right\}^{1/p} \le \frac{\mathsf{L}_{p}}{n} \quad (p=1,2). \label{eq:postn-concentration}
\]
Then
\[
\pwassSimple{p}{\approxdist_{\mathrm{Laplace}}}{\postdist}  \le \frac{\mathsf{L}_{p}\mathsf{M}}{\alpha n}.
\]
\end{proposition}
\cref{eq:postn-concentration} essentially requires that $\postdist_{n}$ concentrates in a ball of radius $O(n^{-1})$ centered at $\mapn$,
which is similar to the conditions appearing elsewhere in the Laplace approximation 
literature such as \citet[Theorem 1]{Kass:1990} or \citet[Theorems 7.108 \& 7.116]{Schervish:1995}. 
These earlier results have stronger regularity conditions on both the density, requiring bounds on the 
first six derivatives of $\logpostn$, and the test function, which must be four times continuously differentiable.
In contrast, we require bounds only on the first three derivatives of $\logpostn$, and 
we can consider all Lipschitz test functions in the $p=1$ case. 
In other ways our results are not comparable with the existing literature. 
An advantage of our approach is that only a single Laplace approximation is needed because the Gaussian approximation
can be applied to any test function. 
For a test function $\phi$, the classical approach to the Laplace approximation is to separately approximate the integrals 
$\int \phi(\param) e^{\loglik(\param)}\priordist(\dee\param) \approx \mcI(\phi)$ and 
$\int e^{\loglik(\param)}\priordist(\dee\param) \approx \mcI(1)$, producing the approximation 
$\int \phi(\param)\postdist(\dee\param) \approx \mcI(\phi)/\mcI(I)$.
Thus, separate integral approximations $\mcI(\phi)$ must be computed for each test function, which could
be computationally expensive if $d$ is large or evaluation of $\loglik(\param)$ is slow (for example if there is a large amount of data).
A benefit of the classical approach is that the results cited above guarantee error of $O(n^{-2})$ whereas \cref{prop:asymptotic-Laplace} 
guarantees error of $O(n^{-1})$. 
The moment generating function approach of \citep{Tierney:1989} does not require computing separate integrals 
for each test function, but it does require the test function to be four times continuously differentiable.

\subsection{Hilbert coresets}
Suppose the data $\data$ are conditionally independent given the parameter 
$\param$, so the log-likelihood decomposes as the sum
$\loglik(\param) = \sum_{j=1}^n \loglik_j(\param)$.
The major cost of posterior inference via Monte Carlo methods in this setting
is the $O(N)$ computation required to evaluate $\loglik(\param)$. 
To reduce this cost, a number of authors have suggested using an approximate log-likelihood
given by $\loglik(w, \param) = \sum_{j=1}^n w_j \loglik_j(\param)$,
where $w = (w_1, \dots, w_n)$ is a set of nonnegative weights~\citep[this or a similar idea appear in][]{Madigan:2002,Feldman:2011b,Zhang:2016,Huggins:2016,Lucic:2018,Campbell:2017,Campbell:2018}.
This log-likelihood approximation induces a \emph{coreset posterior approximation}
 for $\postdist$ given by
\[
\approxdist_{w}(\dee \param) = \frac{e^{\loglik(w,\param)}\priordist(\dee\param)}{\marginallik_w},
\]
where $\marginallik_w$ is the normalizing constant for the approximate posterior given weights $w$.
Computing the coreset likelihood approximation $\loglik(w, \param)$ 
takes $O(\|w\|_0)$ time, where  $\|w\|_0 = \sum_{j=1}^n \ind_{(0, \infty)}(w_j)$
is the number of nonzero weights in $w$. Therefore, the cost 
of inference for $\approxdist_w$ may be significantly reduced 
if $\|w\|_0$ is much smaller than $N$.

The main challenge in building a coreset is finding a sparse set of weights
for which $\approxdist_w$ is still a reasonable approximation to $\postdist$.
\citet{Campbell:2017,Campbell:2018} provide iterative algorithms 
for which $(2,\nu)$-Fisher distance decays \emph{exponentially} in the number of nonzero coreset weights.
They show that there exist constants $0 < \beta < 1$ and $C > 0$ such that for all $K \in \nats$,
after $K$ iterations the output satisfies $\|w\|_0 \leq K$ and $\fd{2, \nu}{\approxdist_w}{\postdist} \leq C\beta^K$.
However, \citet{Campbell:2017,Campbell:2018} do not provide any guidance on the effect of the
weighting distribution $\nu$ aside from suggesting it should be ``close'' to the true posterior $\postdist$, 
nor how to link the guarantee on $\fd{2, \nu}{\approxdist_w}{\postdist}$ back to interpretable statistical guarantees.
\cref{prop:coreset} below resolves both of these issues by combining the guarantees of \citet{Campbell:2017,Campbell:2018} with \cref{cor:Wasserstein-error}.
\begin{proposition}\label{prop:coreset}
Suppose that $\loglik_1, \dots, \loglik_n$, and  $\log\textder{\priordist}{\param}$ are continuously differentiable
and that $-\loglik(w, \cdot) - \log\textder{\priordist}{\param}$ is $\alpha$-strongly convex.
Let $B_{1} = \{1+\dchi{\postdist}{\nu}\}^{1/2}$ 
and $B_{2} =  \infnorm{\textder{\postdist}{\nu}}^{1/2}$.
Then for all $K\in\nats$, after $K$ iterations the output satisfies $\|w\|_0 \leq K$ and 
\[
\pwassSimple{p}{\approxdist_{w}}{\postdist} \leq \alpha^{-1}C B_{p} \beta^K \quad (p=1,2).
\]
\end{proposition}

\section{Discussion}
\label{sec:discussion}

Our results suggest that the $p$-Wasserstein distance is a good choice for measuring posterior approximation accuracy because
it implies bounds on the errors of the estimates for the means, maximum absolute deviations, standard deviations, and covariance. 
For methods that can be viewed as approximating the log-likelihood, the $(p, \nu)$-Fisher distance provides an approach for
computing bounds on the $p$-Wasserstein distance, provided either the posterior or its approximation is strongly convex, as least in the tails. 
For example, promising candidates include variational inference methods~\citep{Blei:2017} and
approximate Bayesian computation~\citep{Marjoram:2003,Marin:2011}, which, as shown in \citet{Karabatsos:2018}, can be viewed as an approximate likelihood method. 
Weakening the tail behavior requirements would be useful in some circumstances such as for the complex likelihoods 
tackled by approximate Bayesian computation.
But such a generalization presents substantial challenges because
heavier tails can more strongly influence the values of Lipschitz functions such as the mean as well as uncertainty measures such as the variance. 
Thus, the price for allowing heavier tail behavior is likely to be much weaker bounds or the introduction of additional restrictive conditions. 
Another aspect of our approach worth careful consideration is that the $(p, \nu)$-Fisher distance is a very strong metric. 
As \cref{prop:Fisher-TV-bound} shows, the Fisher distance bounds the total variance distance, not just the Wasserstein distance. 
This property is useful because it leads to error bounds on credible interval estimates, but it also means that the bounds on
Wasserstein distance are, in general, going to be loose. 
A promising alternative approach would be to instead use a kernel Stein discrepancy, which can be viewed as a kernel-smoothed version 
of the Fisher distance~\citep[Proposition 9]{Gorham:2017}. 
However, choosing an optimal kernel and obtaining tight bounds on means and uncertainty estimates using kernel Stein discrepancies 
remain under-explored questions that would need to be addressed.

\section*{Acknowledgement}
The authors thank Daniel Simpson and Arthur Gretton for valuable discussions and many useful references. 
This research was supported in part by an NSF CAREER Award, an ARO
YIP Award, the Office of Naval Research, and a Sloan Research Fellowship. 
M.\ Kasprzak was supported by an EPSRC studentship.

\appendix
\numberwithin{equation}{section}
\def\subsection*{\section}

\subsection*{Proof of \cref{prop:KL-divergence-working}}

Let $d_{H}^{2}(\estdist, \postdist) = \int\{1 - (\textder{\estdist}{\postdist})^{1/2}\}^{2} \dee\postdist$ denote the squared Hellinger distance.
Without loss of generality assume $\mean{\estdist} = 0$. 
It follows from \citet[Lemma 6.37]{Stuart:2010} that 
\[
(\mean{\estdist} - \mean{\postdist})^{2}
&= \mean{\postdist}^{2} 
\le (\var{\estdist} + \var{\postdist} + \mean{\postdist}^{2}) d_{H}^{2}(\estdist, \postdist) 
\]
and hence, solving for $(\mean{\estdist} - \mean{\postdist})^{2}$, that 
\[
(\mean{\estdist} - \mean{\postdist})^{2} 
&\le (\var{\estdist} + \var{\postdist})d_{H}^{2}(\estdist, \postdist) / \{1 - d_{H}^{2}(\estdist, \postdist)\}.
\]
Since $d_{H}^{2}(\estdist, \postdist) \le  \kl{\estdist}{\postdist}$ and $t \mapsto t/(1-t)$ is monotonically increasing for $t \in [0,1)$, 
the first inequality follows. 
The second inequality follows immediately from Pinsker's inequality and the definition of the total variation distance. 

\subsection*{Proof of \cref{prop:KL-divergence-problems}}

(A) Choose $\estdist$ and $\postdist$ to be Gaussians such that $(\mean{\estdist} - \mean{\postdist})^{2} = \var{\estdist}\{\exp(2\delta) - 1\}$
and $\var{\postdist} = \exp(2\delta)\var{\estdist}$.
We then have that 
\[
\lefteqn{\kl{\estdist}{\postdist}} \\
&= 0{\cdot}5\{\var{\estdist}/\var{\postdist} - 1 + \log(\var{\postdist}/\var{\estdist}) + (\mean{\estdist} - \mean{\postdist})^{2}/\var{\postdist}\} \\
&= 0{\cdot}5[\var{\estdist}/\{\exp(2\delta)\var{\estdist}\} - 1 + \log\{\exp(2\delta)\var{\estdist}/\var{\estdist}\} +  \var{\estdist}\{\exp(2\delta) - 1\}/\{\exp(2\delta)\var{\estdist}\}] \\
&= 0{\cdot}5[\exp(-2\delta) - 1 + \log\{\exp(2\delta)\} +  \{\exp(2\delta) - 1\} \exp(-2\delta)] \\
&= \delta 
\]

(B) Let $\estdist$ be a standard Gaussian and let $\postdist = \mcT_{h}$ be a standard $t$-distribution 
with $h$ degrees of freedom.
For $\Param \dist \estdist$, we have
\[
\lefteqn{\kl{\estdist}{\postdist_{v}}} \\
&= \log[\Gamma(h/2)h^{1/2}/\Gamma\{(h + 1)/2\}]
 - 0{\cdot}5 \log(2e)  + 0{\cdot}5(h + 1)\EE\left\{\log\left(1 + \Param^{2}/h\right)\right\} \\
&\le \log[\Gamma(h/2)h^{1/2}/\Gamma\{(h + 1)/2\}]
 - 0{\cdot}5 \log(2e) + 0{\cdot}5(h + 1)\log\left\{1 + \EE(\Param^{2})/h\right\} \\
&= \log[\Gamma(h/2)h^{1/2}/\Gamma\{(h + 1)/2\}]
 - 0{\cdot}5 \log(2e) + 0{\cdot}5(h + 1)\log\left\{1 + 1/h\right\}.  \label{eq:KL-Guassian-t-bound}
\]
For $h = 2$, \cref{eq:KL-Guassian-t-bound} is equal to $0{\cdot}801345\cdots$.
A tedious but straightforward calculation shows that that $\textder{\kl{\estdist}{\mcT_{h}}}{h} > 0$ at $h = 2$.
Since $\kl{\estdist}{\mcT_{h}}$ is a continuous function of $h$, there exists an $\veps > 0$ such that 
for all $h \in [2,2+\veps)$, $\kl{\estdist}{\mcT_{h}} < 0{\cdot}802$.
Moreover, $\var{\mcT_{h}} \to \infty$ as $h \to 2$ from the right, proving the claim.

\subsection*{Proof of \cref{thm:Wasserstein-moment-bounds}}

We begin by considering the case $d=1$, dropping the component indexes from our notation. 
\begin{theorem} \label{thm:one-dim-Wasserstein-moment-bounds}
Assume $d = 1$. 
If $\pwassSimple{1}{\mainmeas}{\estmeas} \le \veps$, then 
$
|\mean{\mainmeas} - \mean{\estmeas}| \le \veps 
$
and
$
|\MAD{\mainmeas} - \MAD{\estmeas}| \le 2\veps. %
$
If in addition $\mainmeas$ has Lebesgue density bounded by $c < \infty$, then
\[
|I_{\mainmeas,a,b} - I_{\estmeas,a,b}| 
	&\le 2(2 c\,\veps)^{1/2} 
	 \qquad (\infty \le a < b \le \infty). %
\]
On the other hand, if $\pwassSimple{2}{\mainmeas}{\estmeas} \le \veps$, then 
$\pwassSimple{1}{\mainmeas}{\estmeas} \le \veps$,
\[
|\std{\mainmeas} - \std{\estmeas}| \le \frac{1}{2} \left(2^{1/2}+6^{1/2}\right)\veps \le 2\veps, \label{eq:stdev-error-bound}
\]
and
\[
|\var{\mainmeas} - \var{\estmeas}| 
&\le 2^{3/2}\min(\std{\mainmeas}, \std{\estmeas})\veps + (1 + 3 \times 2^{1/2})\veps^{2} \label{eq:var-error-bound} \\
&\le 3\min(\std{\mainmeas}, \std{\estmeas})\veps + 5{\cdot}25 \veps^{2}.
\]
\end{theorem}

The proof of \cref{thm:one-dim-Wasserstein-moment-bounds} is deferred to the next section. 
To generalize to the case of $d > 1$, for a random variable $\Param \dist \eta$ on $\reals^{d}$ with
distribution $\eta$ and any vector $v \in \reals^{d}$, let $\mean{\eta,v} = \EE(v^{\top}\Param)$,
$\var{\eta,v} = \EE\{(v^{\top}\Param -\mean{\eta,v})^{2}\}$, 
$\MAD{\eta,v} = \EE(|v^{\top}\Param - \mean{\eta,v}|)$,
and $I_{\eta,v,a,b} = \EE\{\ind_{[a,b]}(v^{\top}\Param)\}$~($\infty \le a < b \le \infty$). 

\begin{corollary} \label{cor:multidim-Wasserstein-moment-bounds}
Let $v \in \reals^{d}$ satisfy $\twonorm{v} \le 1$. 
If $\pwassSimple{1}{\mainmeas}{\estmeas} \le \veps$ then 
$
|\mean{\mainmeas,v} - \mean{\estmeas,v}| \le \veps
$
and
$
|\MAD{\mainmeas,v} - \MAD{\estmeas,v}| \le 2\veps. \label{eq:v-mean-and-v-MAD-bounds}
$
If in addition $\mainmeas$ has Lebesgue density bounded by $c < \infty$, then
\[
|I_{\mainmeas,v,a,b} - I_{\estmeas,v,a,b}| 
	&\le 2(2 c \veps)^{1/2} 
	 \qquad (\infty \le a < b \le \infty).
\]
On the other hand, if $\pwassSimple{2}{\mainmeas}{\estmeas} \le \veps$, then
\[
|\std{\mainmeas,v} - \std{\estmeas,v}| 
	&\le \frac{1}{2} \left(2^{1/2}+6^{1/2}\right)\veps, \\
|\var{\mainmeas,v} - \var{\estmeas,v}| 
	&\le 2^{3/2}\min(\std{\mainmeas,v}, \std{\estmeas,v})\veps + (1 + 3 \times 2^{1/2})\veps^{2}.
\]
\end{corollary}
\bprf
Let $\Param \dist \mainmeas$, let $\Param_{v} = v^{\top}\Param_{(j)}$ and let $\mainmeas_{v}$ denote the 
distribution of $\Param_{v}$. 
Define $\hat{\Param}$, $\hat{\Param}_{v}$, and $\estmeas_{v}$ analogously in terms of $\estmeas$. 
By the Cauchy-Schwarz inequality and the assumption that $\twonorm{v} \le 1$, 
\[
\EE(|\Param_{v} - \hat{\Param}_{v}|^{p})
= \EE(|v^{\top}\Param - v^{\top}\hat{\Param}|^{p})
\le \EE(\statictwonorm{\Param - \hat{\Param}}^{p}).
\]
Hence $\pwassSimple{p}{\mainmeas_{v}}{\estmeas_{v}} \le \pwassSimple{p}{\mainmeas}{\estmeas}$. 
The corollary now follows from \cref{thm:one-dim-Wasserstein-moment-bounds}. 
\eprf
\begin{lemma} \label{lem:sup-v-equalities}
For probability measures $\xi, \mainmeas, \estmeas$, we have  
$\twonorm{\mean{\mainmeas} - \mean{\estmeas}} = \sup_{\twonorm{v} \le 1}|\mean{\mainmeas,v} - \mean{\estmeas,v}|$, 
$\twonorm{\Sigma_{\xi}} = \sup_{\twonorm{v} \le 1} \var{\xi,v}$, and  
$\twonorm{\Sigma_{\mainmeas} - \Sigma_{\estmeas}} = \sup_{\twonorm{v} \le 1} |\var{\mainmeas,v} - \var{\estmeas,v}|$. 
\end{lemma}
\bprf
The first result follows since $\mean{\mainmeas,v} - \mean{\estmeas,v} = v^{\top}(\mean{\mainmeas} - \mean{\estmeas})$
and for any $w \in \reals^{d}$, $\sup_{\twonorm{v} \le 1} v^{\top}w = \twonorm{w}$. 
For the second result, since $\Sigma_{\xi}$ is positive semi-definite,
\[
\twonorm{\Sigma_{\xi}}
&= \sup_{\twonorm{v} \le 1} v^{\top}\Sigma_{\xi}v  
= \sup_{\twonorm{v} \le 1}\EE\{v^{\top}(\Xvec - \meanvec{\xi})(\Xvec - \meanvec{\xi})^{\top}v\} 
= \sup_{\twonorm{v} \le 1}\var{\xi,v};
\]
The third result follows by an analogous argument.
\eprf

By taking $v = e_{i}$, the $i$th canonical basis vector of $\reals^{d}$, 
\cref{cor:multidim-Wasserstein-moment-bounds} implies the bounds in \cref{thm:Wasserstein-moment-bounds}
on $|\MAD{\mainmeas,i} - \MAD{\estmeas,i}|$, $|I_{\mainmeas,i,a,b} - I_{\estmeas,i,a,b}|$, and 
$|\std{\mainmeas,i} - \std{\estmeas,i}|$.
\cref{cor:multidim-Wasserstein-moment-bounds,lem:sup-v-equalities} yield the bounds in \cref{thm:Wasserstein-moment-bounds}
on $\twonorm{\mean{\mainmeas} - \mean{\estmeas}}$ and $\twonorm{\Sigma_{\mainmeas} - \Sigma_{\estmeas}}$.

\subsection*{Proof of \cref{thm:one-dim-Wasserstein-moment-bounds}}

Throughout we will always assume that $\Param \dist \mainmeas$ and $\hat{\Param} \dist \estmeas$
are distributed according to the optimal coupling for the $p$-Wasserstein distance under consideration.
We will also assume without loss of generality that $\mean{\mainmeas} = 0$ since if not
we could consider the random variables $\Param' = \Param - \mean{\mainmeas}$ and $\hat{\Param}' = \hat{\Param} - \mean{\mainmeas}$
instead. 

Assume $\pwassSimple{1}{\mainmeas}{\estmeas} \le \veps$.
By \cref{eq:1-Wasserstein-dual}, for any Lipschitz function $\phi$,
\[
|\EE(\phi(\Param) - \phi(\hat{\Param}))| \le \veps\norm{\phi}_{L}. 
\]
Hence, taking $\phi(t) = t$, we have that 
$
|\mean{\mainmeas} - \mean{\estmeas}| = |\mean{\estmeas}| \le \veps. 
$
For the mean absolute deviation, using the fact that $\phi(t) = |t|$
is 1-Lipschitz, we have 
\[
|\MAD{\mainmeas} - \MAD{\estmeas}|
&= |\EE(|\Param| - |\hat{\Param} - \mean{\estmeas}|)| 
\le |\EE(|\Param| - |\hat{\Param}|)| + |\mean{\estmeas}| 
\le 2\veps.
\]
\cref{eq:credible-interval-error-bound} follows immediately from the 1-Wasserstein distance
bound on the Kolmogorov distance~\citep[Appendix C]{Nourdin:2012}.

Assume $\pwassSimple{2}{\mainmeas}{\estmeas} \le \veps$.
By Jensen's inequality 
$\pwassSimple{1}{\mainmeas}{\estmeas} \le \veps$ as well. 
Let $\varsigma_{\mainmeas}^{2} = \EE(\Param^{2}) = \var{\mainmeas}$ and $\varsigma_{\estmeas}^{2} = \EE(\hat{\Param}^{2})$.
It follows from the Cauchy-Schwarz inequality that 
\[
|\varsigma_{\mainmeas}^{2} - \varsigma_{\estmeas}^{2}|
&= \big|\EE\big(\Param^{2} - \hat{\Param}^{2}\big)\big| 
= \big|\EE\left\{\left(\Param - \hat{\Param}\right)\left(\Param + \hat{\Param}\right)\right\}\big| \\
&\le \EE\left\{\left(\Param - \hat{\Param}\right)^{2}\right\}^{1/2}\EE\left\{\left(\Param + \hat{\Param}\right)^{2}\right\}^{1/2}
\le 2^{1/2}\veps \EE\big(\Param^{2} + \hat{\Param}^{2}\big)^{1/2} \\
&\le 2^{1/2}\veps (\varsigma_{\mainmeas} + \varsigma_{\estmeas}).
\]
Since $|\varsigma_{\mainmeas}^{2} - \varsigma_{\estmeas}^{2}| = |\varsigma_{\mainmeas} - \varsigma_{\estmeas}| (\varsigma_{\mainmeas} + \varsigma_{\estmeas})$, 
it follows that 
\[
|\varsigma_{\mainmeas} - \varsigma_{\estmeas}| &\le 2^{1/2}\veps.  \label{eq:bar-sigma-diff}
\]
Using \cref{eq:bar-sigma-diff}, we also have
\[
|\var{\mainmeas} - \var{\estmeas}|
    &= |\varsigma_{\mainmeas}^{2} - \varsigma_{\estmeas}^{2} + \mean{\estmeas}^{2}| 
    \le |\varsigma_{\mainmeas}^{2} - \varsigma_{\estmeas}^{2}| + |\mean{\estmeas}^{2}|
    \le 2^{1/2}\veps (\varsigma_{\mainmeas} + \varsigma_{\estmeas}) + \veps^{2} \label{eq:var-error-preliminary} \\
|\std{\mainmeas} - \std{\estmeas}| 
    &\le 2^{1/2}\veps + \frac{\veps^{2}}{\varsigma_{\mainmeas} + \varsigma_{\estmeas}}. \label{eq:stdev-error-preliminary}
\]
If $\max(\std{\mainmeas}, \std{\estmeas}) \le \frac{1}{2} \left(2^{1/2}+6^{1/2}\right) \veps$, then 
clearly $|\std{\mainmeas} - \std{\estmeas}| \le  \frac{1}{2} \left(2^{1/2}+6^{1/2}\right) \veps$. 
Otherwise $\frac{\veps^{2}}{\varsigma_{\mainmeas} + \varsigma_{\estmeas}} \le \frac{2\veps}{2^{1/2}+6^{1/2}}$  
and so, using \cref{eq:stdev-error-preliminary}, we have 
\[
|\std{\mainmeas} - \std{\estmeas}| 
&\le 2^{1/2}\veps + \frac{2\veps}{2^{1/2}+6^{1/2}}
= \frac{1}{2} \left(2^{1/2}+6^{1/2}\right) \veps.
\]
Hence we conclude unconditionally that 
$
|\std{\mainmeas} - \std{\estmeas}| 
\le  \frac{1}{2} \left(2^{1/2}+6^{1/2}\right) \veps.
$ 
Starting with \cref{eq:var-error-preliminary} and using \cref{eq:bar-sigma-diff}, we have
\[
|\var{\mainmeas} - \var{\estmeas}| 
&\le 2^{1/2}\veps (\varsigma_{\mainmeas} + \varsigma_{\estmeas}) + \veps^{2}  
= 2^{1/2}\veps\left \{\std{\mainmeas} + (\var{\estmeas} + \mean{\estmeas}^{2})^{1/2}\right\} + \veps^{2}  \\
&\le 2^{1/2}\veps (2\std{\mainmeas} + 3\veps) + \veps^{2} 
= 2^{3/2}\,\std{\mainmeas}\veps + (1 + 3 \times 2^{1/2})\veps^{2}.
\]

\subsection*{Proof of \cref{thm:p-Wasserstein-sde-error} ($p=1$ case)}

Let $\drift = \grad\mainPE$ and $\hat{\drift} = \grad\estPE$. 
Consider the following $\reals^{d}$-valued diffusions with respect to a $d$-dimensional Wiener process $\Wiener{}$
\[
\dXt{t} = \drift(\Xt{t})\dt + 2^{1/2}\dWt{t},\quad
\dYt{t} = \hat{\drift}(\Yt{t})\dt + 2^{1/2}\dWt{t},
\]
which have unique stationary measures $\mainmeas$ and $\estmeas$, respectively. We couple them using the ``coupling of marching soldiers''~\citep[Example 2.16]{Chen:2005}
\[ 
\dee(\Xt{t},\Yt{t})= \{\drift(\Xt{t}),\hat{\drift}(\Yt{t})\}\dt + 2^{1/2}\dee(\Wiener{t},\Wiener{t}),\label{sde}
\]
and assume that the processes $\Xt{}$ and $\Yt{}$ are both started at stationarity (with $\Xt{0}\sim\mainmeas$ and $\Yt{0}\sim\estmeas$).
It follows from the $\alpha$-strong convexity of $\estPE$ that $\hat{\drift}$ satisfies
\[
\inner{\hat{\drift}(\param)-\hat{\drift}(\param')}{\param-\param'}\leq -\alpha\twonorm{\param-\param'}^2 \quad (\theta, \theta' \in \reals^{d}). \label{eq:one-sided-Lipschitz}
\]
Using \cref{eq:one-sided-Lipschitz} and the Cauchy-Schwarz inequality, we have that 
\[
\lefteqn{\inner{\drift(\Xt{t}) - \hat{\drift}(\Yt{t})}{\Xt{t} - \Yt{t}}}  \\
&= \inner{\hat{\drift}(\Xt{t}) - \hat{\drift}(\Yt{t})}{\Xt{t} - \Yt{t}} + \inner{\drift(\Xt{t}) - \hat{\drift}(\Xt{t})}{\Xt{t} - \Yt{t}} \\
&\le -\alpha\twonorm{\Xt{t} - \Yt{t}}^{2} + \staticnorm{\drift(\Xt{t}) - \hat{\drift}(\Xt{t})}_{2}\twonorm{\Xt{t} - \Yt{t}}. \label{eq:drift-difference-inequality}
\]
In order to obtain the estimate for $\pwassSimple{1}{\mainmeas}{\estmeas}$, we will follow a strategy similar to the one used in a proof of Tanaka's formula~\citep[see][Exercise 4.10]{Oksendal:2003}. For any $\epsilon>0$, let us consider $g_{\epsilon}:\reals^d\to\reals$ given by:
$$g_{\epsilon}(x)=\begin{cases}
\twonorm{x},\quad &\twonorm{x}\geq \epsilon\\
\frac{1}{2}\left(\epsilon+\frac{\twonorm{x}^2}{\epsilon}\right),\quad &\twonorm{x}<\epsilon.
\end{cases}$$
Applying It\^{o}'s formula~\citep[Theorem 4.2.1]{Oksendal:2003} to the SDE (\ref{sde}) and function ${g:[0,1]\times\reals^d\times\reals^d\to\reals}$, given by $g(t,x,y)=e^{\alpha t}g_{\epsilon}(x-y)$, we obtain
\[
\lefteqn{e^{\alpha t}g_{\epsilon}(\Xt{t}-\Yt{t})}\\
&=g_{\epsilon}(\Xt{0}-\Yt{0})+\int_0^te^{\alpha s}2^{1/2}\inner{\nabla g_{\epsilon}(\Xt{s}-\Yt{s})}{\dWt{s}-\dWt{s}}+\int_0^t\alpha e^{\alpha s}g_{\epsilon}(\Xt{s}-\Yt{s})\dee s\nonumber\\
&\phantom{=~}+\int_0^te^{\alpha s}\binner{\frac{\Xt{s}-\Yt{s}}{\twonorm{\Xt{s}-\Yt{s}}}}{\drift(\Xt{s})-\hat{\drift}(\Yt{s})}\mathbb{I}[\twonorm{\Xt{s}-\Yt{s}}\geq \epsilon]\dee s\nonumber\\
&\phantom{=~}+\int_0^t e^{\alpha s}\binner{\frac{\Xt{s}-\Yt{s}}{\epsilon}}{\drift(\Xt{s})-\hat{\drift}(\Yt{s})}\mathbb{I}[\twonorm{\Xt{s}-\Yt{s}}<\epsilon]\dee s\nonumber\\
&\phantom{=~}+4\sum_{i=1}^{d}(1-1)\int_0^te^{\alpha s}\left(\frac{1}{\twonorm{\Xt{s}-\Yt{s}}}-\frac{\left(\Xt{s}^{(i)}-\Yt{s}^{(i)}\right)^2}{\twonorm{\Xt{s}-\Yt{s}}^3}\right)\mathbb{I}[\twonorm{\Xt{s}-\Yt{s}}\geq\epsilon]\dee s\\
&\phantom{=~}+2\sum_{i=1}^{2d}(1-1)\int_0^t\frac{e^{\alpha s}}{\epsilon}\mathbb{I}[\twonorm{\Xt{s}-\Yt{s}}<\epsilon]\dee s\\
&=g_{\epsilon}(\Xt{0}-\Yt{0})+\int_0^t\alpha e^{\alpha s}\twonorm{\Xt{s}-\Yt{s}}\mathbb{I}[\twonorm{\Xt{s}-\Yt{s}}\geq\epsilon]\dee s\\
&\phantom{=~}+\frac{1}{2}\int_0^t\alpha e^{\alpha s}\left(\epsilon+\frac{\twonorm{\Xt{s}-\Yt{s}}^2}{\epsilon}\right)\mathbb{I}[\twonorm{\Xt{s}-\Yt{s}}<\epsilon]\dee s\\
&\phantom{=~}+\int_0^te^{\alpha s}\binner{\frac{\Xt{s}-\Yt{s}}{\twonorm{\Xt{s}-\Yt{s}}}}{\drift(\Xt{s})-\hat{\drift}(\Yt{s})}\mathbb{I}[\twonorm{\Xt{s}-\Yt{s}}\geq \epsilon]\dee s\\
&\phantom{=~}+\int_0^t e^{\alpha s}\binner{\frac{\Xt{s}-\Yt{s}}{\epsilon}}{\drift(\Xt{s})-\hat{\drift}(\Yt{s})}\mathbb{I}[\twonorm{\Xt{s}-\Yt{s}}<\epsilon]\dee s.
\]
Using Eq. (\ref{eq:drift-difference-inequality}), we obtain
\[
\lefteqn{e^{\alpha t}g_{\epsilon}(\Xt{t}-\Yt{t})}\nonumber\\
&\leq g_{\epsilon}(\Xt{0}-\Yt{0})+\int_0^te^{\alpha s}\staticnorm{\drift(\Xt{s})-\hat{\drift}(\Xt{s})}_{2}\mathbb{I}[\twonorm{\Xt{s}-\Yt{s}}\geq\epsilon]\dee s\\
&\phantom{=~}+\frac{1}{2}\int_0^t\alpha e^{\alpha s}\left(\epsilon+\frac{\epsilon^2}{\epsilon}\right)\mathbb{I}[\twonorm{\Xt{s}-\Yt{s}}<\epsilon]\dee s\\
&\phantom{=~}+\int_0^t e^{\alpha s}\left(-\alpha\frac{\twonorm{\Xt{s}-\Yt{s}}^2}{\epsilon}+\frac{\staticnorm{\drift(\Xt{s})-\hat{\drift}(\Xt{s})}_{2}\twonorm{X_s-Y_s}}{\epsilon}\right)\mathbb{I}[\twonorm{\Xt{s}-\Yt{s}}<\epsilon] \dee s\\
&\leq g_{\epsilon}(\Xt{0}-\Yt{0})+\int_0^te^{\alpha s}\staticnorm{\drift(\Xt{s})-\hat{\drift}(\Xt{s})}_{2}\mathbb{I}[\twonorm{\Xt{s}-\Yt{s}}\geq\epsilon]\dee s+\epsilon(e^{\alpha t}-1)\\
&\phantom{=~}-\epsilon(e^{\alpha t}-1)+\int_0^te^{\alpha s}\staticnorm{\drift(\Xt{s})-\hat{\drift}(\Xt{s})}_{2}\mathbb{I}[\twonorm{X_s-Y_s}<\epsilon]\dee s\\
&= g_{\epsilon}(\Xt{0}-\Yt{0})+\int_0^te^{\alpha s}\staticnorm{\drift(\Xt{s})-\hat{\drift}(\Xt{s})}_{2}\dee s.
\]
Taking $\epsilon\to 0$ and taking expectations on both sides (at a fixed time $t$, with respect to everything that is random), we obtain
\[
e^{\alpha t}\EE(\twonorm{\Xt{t}-\Yt{t}}) 
\leq \EE(\twonorm{\Xt{0}-\Yt{0}}) + \alpha^{-1}(e^{\alpha t}-1)\EE(\staticnorm{\drift(\Xt{t})-\hat{\drift}(\Xt{t})}_{2}),
\]
which follows on the assumption that the process $t\mapsto \Xt{t}$ was started at stationarity.
Dividing by $e^{\alpha t}$, taking $t\to\infty$ and retaining the assumption that both processes $t\mapsto \Yt{t}$ and $t\mapsto \Xt{t}$ are at stationarity, we obtain
\[
\pwassSimple{1}{\mainmeas}{\estmeas}\leq \alpha^{-1}\fd{1,\mainmeas}{\mainmeas}{\estmeas}.
\]

\subsection*{Proofs of \cref{thm:p-Wasserstein-sde-error} ($p=2$ case) and \cref{thm:p-Wasserstein-sde-error-general}}

In the notation of \citet{Bolley:2012}, the $p=2$ case of \cref{eq:p-Wasserstein-sde-error} would be written
\[
\pwassSimple{2}{\mainmeas}{\estmeas} \le C^{-1} I(\mainmeas|\estmeas)^{1/2} \label{eq:ineq-2}
\]
where $C = 1/\alpha$ and $I(\mainmeas|\estmeas) =  \fd{\mainmeas}{\mainmeas}{\estmeas}^{2}$. 
When \cref{eq:ineq-2} holds for all absolutely continuous $\mainmeas$, \citet{Bolley:2012} say that $\estmeas$ satisfies a $WI(C)$ inequality (see p.~2450). 
The definition of  $I(\mainmeas|\estmeas)$, the Fisher information of $\mainmeas$ with respect to $\estmeas$, is~\citep[p.~2436]{Bolley:2012}
\[
I(\mainmeas|\estmeas) = \int \frac{\twonorm{\grad h}^{2}}{h} \dee \estmeas,
\]
where $h = \dee \mainmeas/\dee \estmeas$. 
We abuse notation and use $\mainmeas$ and $\estmeas$ to denote their respective densities. 
Then we verify the equivalence of \cref{eq:p-Wasserstein-sde-error} ($p=2$) and \cref{eq:ineq-2}:
\[
I(\mainmeas|\estmeas) 
&= \int \frac{\twonorm{\grad (\mainmeas/\estmeas)}^{2}}{\mainmeas/\estmeas} \dee \estmeas 
= \int \frac{\estmeas^{2}}{\mainmeas^{2}}\twonorm{(1/\mainmeas)\grad \estmeas - (\mainmeas/\estmeas^{2}) \grad \estmeas}^{2} \dee \mainmeas \\
&= \int \twonorm{(\grad \mainmeas)/\mainmeas - (\grad \estmeas)/\estmeas}^{2} \dee \mainmeas 
= \int \twonorm{\grad \log\mainmeas - \grad \log \estmeas}^{2} \dee \mainmeas \\
&= \fd{\mainmeas}{\mainmeas}{\estmeas}^{2}.
\]
\citet{Bolley:2012} also introduce what they call the $WJ(C)$ inequality (the details of what it is are not
important for our purpose). 
\citet[Prop.~3.10]{Bolley:2012} states that $WJ(C) \implies WI(C)$, so in order to verify \cref{eq:p-Wasserstein-sde-error} ($p=2$) it
suffices to show that $\estmeas$ satisfies the $WJ(C)$ inequality condition. 
\citet[Lem.~3.3]{Bolley:2012} states that if $\estPE$ is $C$-strongly convex then $\estmeas$ satisfies the $WJ(C)$ inequality condition,
which concludes the proof of \cref{thm:p-Wasserstein-sde-error} ($p=2$ case).
\citet[Prop.~3.4]{Bolley:2012} states that if $\estPE$ satisfies the hypotheses in \cref{thm:p-Wasserstein-sde-error-general}
then $\estmeas$ satisfies a $WJ$ inequality with constant $C = C(V, R, K)$, which concludes the proof of \cref{thm:p-Wasserstein-sde-error-general}.

\subsection*{Proof of \cref{cor:Wasserstein-error}}

Using the Cauchy-Schwarz inequality, we have
\[
\fd{1,\nu}{\mainmeas}{\estmeas}^{2}
&= \left(\int \staticnorm{\drift(\param) - \hat{\drift}(\param)} \mainmeas(\dee\param)\right)^{2} 
= \left(\int \staticnorm{\drift(\param) - \hat{\drift}(\param)} \der{\mainmeas}{\nu}(\param) \nu(\dee\param)\right)^{2} \\
&\le \int \staticnorm{\drift(\param) - \hat{\drift}(\param)}^{2}  \nu(\dee\param)  \int \der{\mainmeas}{\nu}(\param)^{2} \nu(\dee\param) 
= \fd{2,\nu}{\mainmeas}{\estmeas}^{2} \int \der{\mainmeas}{\nu}(\param) \mainmeas(\dee\param) \\
&= \fd{2,\nu}{\mainmeas}{\estmeas}^{2} (1 + \dchi{\mainmeas}{\nu}).  \label{eq:1-nu-F-dist-bound}
\]
\cref{eq:p-Wasserstein-error} ($p=1$) follows by combining \cref{eq:p-Wasserstein-sde-error} ($p = 1$) and \cref{eq:1-nu-F-dist-bound}.
Now using \Holder's inequality, we have
\[
\fd{2,\nu}{\mainmeas}{\estmeas}^{2}
&= \int \staticnorm{\drift(\param) - \hat{\drift}(\param)}^{2} \mainmeas(\dee\param) 
= \int \staticnorm{\drift(\param) - \hat{\drift}(\param)}^{2} \der{\mainmeas}{\nu}(\param) \nu(\dee\param) \\
&\le \int \staticnorm{\drift(\param) - \hat{\drift}(\param)}^{2} \nu(\dee\param)  \infnorm{\der{\mainmeas}{\nu}}. \label{eq:2-nu-F-dist-bound}
\]
\cref{eq:p-Wasserstein-error} ($p=2$) follows by  combining \cref{eq:p-Wasserstein-sde-error} ($p = 2$) and \cref{eq:2-nu-F-dist-bound}.

\subsection*{Proofs of Propositions \ref{prop:normal-bound-tightness} and \ref{prop:t-dist-bound-tightness}}

For \cref{prop:normal-bound-tightness}, it is easy to check that $\infnorm{\textder{\mainmeas}{\nu}} = C$. 
A straightforward but tedious calculation shows that 
\[
\fd{2,\nu}{\mainmeas}{\estmeas} = \frac{\rho (\std{\mainmeas}^3-\std{\estmeas}^2 \std{\mainmeas} )^2+[\std{\mainmeas} ^2 (\Delta\mean{} -\epsilon )+\std{\estmeas}^2 \epsilon]^2}{\std{\estmeas}^4 \std{\mainmeas} ^4}.
\]
Then, using the fact that $\estPE$ is $\std{\estmeas}^{-2}$-strongly convex, the result follows after some further algebra. 

\cref{prop:t-dist-bound-tightness} follows by similar arguments. 

\subsection*{Proof of \cref{prop:Fisher-TV-bound}}

\citet[Prop.~3.9]{Bolley:2012} states that if $\estmeas$ satisfies the $WJ(C)$ inequality condition then 
it satisfies a Poincar\'e inequality with the constant $C^{-1}$: for all measurable $\phi : \reals^{d} \to \reals$ 
such that $\estmeas(\phi) = 0$, 
\[
\int \phi(\param)^{2} \estmeas(\dee \param) 
\le C^{-1} \int \twonorm{\grad \phi(\param)}^{2} \estmeas(\dee \param),
\]
Following the approach of \citet[Lem.~E.2]{Johnson:2004a}, we can show that if $\estmeas$ satisfies 
a Poincar\'e inequality then the Hellinger distance $d_{H}(\mainmeas, \estmeas) = \{\int|\mainmeas(\param)^{1/2} - \estmeas(\param)^{1/2}|^{2}\dee \param\}^{1/2}$
can be bounded by the Fisher distance. 
\begin{lemma} \label{lem:Poincare}
If $\estmeas$ satisfies a Poincar\'e inequality with constant $C^{-1}$ then for all absolutely continuous $\mainmeas$,
\[
d_{H}(\mainmeas, \estmeas) \le (2C)^{-1/2}\fd{2,\mainmeas}{\mainmeas}{\estmeas}.
\]
\end{lemma}
\cref{prop:Fisher-TV-bound} follows from \cref{lem:Poincare}, \cref{eq:2-nu-F-dist-bound} and
the fact that $d_{\textrm{TV}}(\mainmeas, \estmeas) \le 2 d_{H}(\mainmeas, \estmeas)$.

\subsection*{Proof of \cref{lem:Poincare}}
 
Applying the Poincar\'e inequality to the function $g(\param) = \{\mainmeas(\param)/\estmeas(\param)\}^{1/2}$ we have 
\[
\int (g(\param) - m)^{2} \estmeas(\dee \param) 
\le C^{-1} \int \twonorm{\grad g(\param)}^{2} \estmeas(\dee \param),
\]
where $m = \int g(\param) \estmeas(\dee \param) = \int \sqrt{\mainmeas(\param)\estmeas(\param)}\dee \param$. 
The left hand side can be written as 
\[
\int (g(\param) - m)^{2} \estmeas(\dee \param) 
&= \int g(\param)^{2} \estmeas(\dee \param) - m^{2}
= 1 - m^{2}.
\]
while for the right hand side we have
\[
\int \twonorm{\grad g(\param)}^{2} \estmeas(\dee \param) 
&= \int \frac{1}{4}\der{\estmeas}{\mainmeas}(\param)\twonorm{\grad(\mainmeas/\estmeas)(\param)}^{2} \estmeas(\dee \param) 
= \frac{1}{4}\fd{2,\mainmeas}{\mainmeas}{\estmeas}^{2}.
\]
Next we rewrite the squared Hellinger distance as 
\[
d_{H}(\mainmeas, \estmeas)^{2}
&= \int|\sqrt{\mainmeas(\param)} - \sqrt{\estmeas(\param)}|^{2}\dee \param \\
&= \int(\mainmeas(\param)  + \estmeas(\param) - 2 \sqrt{\mainmeas(\param)\estmeas(\param)})\dee \param \\
&= 2(1 - m).
\]
Since $m \in [0,1]$, conclude that 
\[
d_{H}(\mainmeas, \estmeas)^{2}
= 2(1 - m)
\le 2(1 - m^{2})
\le (2C)^{-1} \fd{2,\mainmeas}{\mainmeas}{\estmeas}^{2}.
\]

\subsection*{Proof of \cref{prop:nonasymptotic-Laplace}}

Let $\drift(\param) = \grad \log \postdensity(\param)$ %
and $\drift_{\mathrm{Laplace}}(\param) = -\maphessian(\param - \map)$. 
By Taylor's theorem, the $i$th component of $\drift(\param)$ can be rewritten as 
\[
\drift_{i}(\param) 
&= \partial_{i} \log \postdensity(\map) + \grad \partial_{i}\log \postdensity(\map)^{\top}(x - \map) + R(\partial_{i} \log \postdensity, \param) \\
&=  \grad \partial_{i}\log \postdensity(\map)^{\top}(x - \map) + R(\partial_{i} \log \postdensity, \param), \label{eq:drift-taylor-expansion}
\]
where 
\[
R(\phi, \param) = (\param - \map)^{\top}\left\{\int_{0}^{1}(1-t)(\hessian \phi)(\map + t(\param - \map))\,\dee t\right\}(\param - \map).
\]
Hence,
\[
\twonorm{\drift(\param) - \drift_{\mathrm{Laplace}}(\param)}^{2}
&= \sum_{i=1}^{d} R(\partial_{i} \log \postdensity, \param)^{2} \\
&\le \sup_{t \in [0,1]}\sum_{i=1}^{d}\twonorm{\param - \map}^{4}\twonorm{\{\hessian\partial_{j}\log\postdensity\}\{\map + t(\param - \map)\}}^{2} \\
&\le M^{2}\twonorm{\param - \map}^{4}. \label{eq:Laplace-drift-difference-bound}
\]
Let $\Lambda$ denote the diagonal matrix with $\Lambda_{ii} = \lambda_{i}~(i=1,\dots,d)$,
$X \dist \distNorm(0, \Lambda)$, and $\Param \dist \approxdist_{\mathrm{Laplace}}$.
Then $\twonorm{X}$ is equal in distribution to $\twonorm{\Param - \map}$. 
It is straightforward to compute expected powers of the norm of $X$.
\begin{lemma} \label{lem:Gaussian-norms}
For $X$ defined above, $\EE(\twonorm{X}^{2}) = \onenorm{\lambda}$
and $\EE(\twonorm{X}^{4}) = 2 \twonorm{\lambda}^{2} + \onenorm{\lambda}^{2}$.
\end{lemma}
Using \cref{lem:Gaussian-norms},
the result follows from \cref{thm:p-Wasserstein-sde-error} with $\mainmeas = \approxdist_{\mathrm{Laplace}}$
and $\estmeas = \postdist$.

\subsection*{Proof of \cref{prop:asymptotic-Laplace}}

The proof is essentially identical to that of \cref{prop:nonasymptotic-Laplace}. 
However, we apply \cref{thm:p-Wasserstein-sde-error} with $\mainmeas = \postdist$
and $\estmeas = \approxdist_{\mathrm{Laplace}}$.
By assumption $-\log \approxdensity_{\mathrm{Laplace}}$ is $\alpha$-strongly convex 
and for $\Param_{n} \dist \postdist_{n}$, by Eqs.~\eqref{eq:postn-concentration} and \eqref{eq:Laplace-drift-difference-bound},
\[
\EE[\twonorm{\drift(\param) - \drift_{\mathrm{Laplace}}(\param)}^{p}]^{1/p} \le  M L_{p}/n,
\]
proving the result. 

\subsection*{Proof of \cref{lem:Gaussian-norms}}

We use the fact that $\EE(X_{i}^{2}) = \lambda_{i}$, $\EE(X_{i}^{4}) = 3\lambda_{i}^{2}$,
and $\EE(X_{i}^{2k-1}) = 0~(i=1,\dots,d; k \in \nats)$. 
For the first equality, we have
\[
\EE(\twonorm{X}^{2}) 
&= \EE\left(\textstyle\sum_{i=1}^{d}X_{i}^{2}\right)
= \textstyle\sum_{i=1}^{d} \lambda_{i} 
= \onenorm{\lambda}.  \label{eq:second-moment}
\]
For the second equality, we have
\[
\EE(\twonorm{X}^{4}) 
&= \EE\left\{\left(\textstyle\sum_{i=1}^{d}X_{i}^{2}\right)^{2}\right\}
= \EE\left(\textstyle\sum_{i=1}^{d}X_{i}^{4} + 2\textstyle\sum_{i=1}^{d}\sum_{j=1}^{i-1}X_{i}^{2}X_{j}^{2}\right) \\
&= 3\textstyle\sum_{i=1}^{d}\lambda_{i}^{2} + 2\textstyle\sum_{i=1}^{d}\sum_{j=1}^{i-1}\lambda_{i}\lambda_{j} 
= 2\textstyle\sum_{i=1}^{d}\lambda_{i}^{2} + \textstyle\sum_{i=1}^{d}\sum_{j=1}^{d}\lambda_{i}\lambda_{j} \\
&= 2 \twonorm{\lambda}^{2} + \onenorm{\lambda}^{2}.  \label{eq:fourth-moment}
\]

\opt{extra-Laplace}{
\section*{Additional Laplace approximation result}

\TBD{This section can probably be removed for submission}

\begin{proposition} \label{prop:Laplace-approximation}
Assume that $-\log\postdensity$ is four times continuously differentiable and $\alpha$-strongly convex
and that $L_{i} = \sup_{\param \in \reals^{d}} \twonorm{\derivop{3}\partial_{i}\log\postdensity[\param]} < \infty~(i=1,\dots,d)$.
Let $\lambda$ denote the eigenvalues of $(\maphessian)^{-1}$, $C_{1}(\lambda) = \onenorm{\lambda}$, 
$C_{2}(\lambda) = (2 \twonorm{\lambda}^{2} + \onenorm{\lambda}^{2})^{1/2}$, 
and $C_{3}(\lambda) = (8 \norm{\lambda}_{3}^{3} + 6\twonorm{\lambda}^{2}\onenorm{\lambda} + \onenorm{\lambda}^{3})^{1/3}$.
Then, letting $M^{\star}_{i} = \twonorm{(\hessian\partial_{i}\log\postdensity)(\map)}$,
\[
\pwassSimple{1}{\approxdist_{\mathrm{Laplace}}}{\postdist}  &\le \alpha^{-1}\left\{\frac{1}{2}\onenorm{M^{\star}}C_{1}(\lambda) + \frac{1}{6}\onenorm{L}C_{2}(\lambda)^{3/2}\right\}  \\
\pwassSimple{2}{\approxdist_{\mathrm{Laplace}}}{\postdist}  &\le \alpha^{-1}\left[\frac{1}{2}\twonorm{M^{\star}}C_{2}(\lambda) + \frac{1}{6^{1/2}}\{(M^{\star})^{\top}L\}^{1/2} C_{3}(\lambda)^{5/4} + \frac{1}{6}\twonorm{L}C_{3}(\lambda)^{3/2}\right]
\]
\end{proposition}

For a 3-tensor $T$, define the operator norm $\opnorm{T} = \sup_{\param \in \reals^{d} \st \twonorm{\param}=1}\twonorm{T[\param]}$.
Let $\drift(\param) = \grad \log \postdensity(\param) = \grad \loglik(\param) + \grad \log \priordensity(\param)$,
$\drift_{\mathrm{Laplace}}(\param) = -\maphessian(\param - \map)$, and $\delta = \param - \map$. 
By Taylor's theorem, the $i$th component of $\drift(\param)$ can be rewritten as 
\[
\drift_{i}(\param) 
&= \partial_{i} \log \postdensity(\map) + \grad \partial_{i}\log \postdensity(\map)^{\top}\delta + \frac{1}{2}\delta^{\top}(\hessian\partial_{i}\log\postdensity)(\map)\delta  + R(\partial_{i}\log\postdensity, \param) \\
&=  \grad \partial_{i}\log \postdensity(\map)^{\top}(\param - \map) + \frac{1}{2}\delta^{\top}(\hessian\partial_{i}\log\postdensity)(\map)\delta + R(\partial_{i} \log \postdensity, \param), \label{eq:drift-taylor-expansion}
\]
where the remainder term satisfies
\[
R(\phi, \param) 
&= \frac{1}{2}\sum_{i=1}^{d}\sum_{j=1}^{d}\sum_{k=1}^{d}\delta_{i}\delta_{j}\delta_{k} \int_{0}^{1}(1-t)^{2}\partial_{i}\partial_{j}\partial_{k}\phi(\map + t\delta)\dee t \\
&\le  \frac{1}{6}\sup_{\param' \in \reals^{d}}\sum_{i=1}^{d}\sum_{j=1}^{d}\sum_{k=1}^{d}\delta_{i}\delta_{j}\delta_{k} \partial_{i}\partial_{j}\partial_{k}\phi(\param') 
=  \frac{1}{6}\sup_{\param' \in \reals^{d}}(\derivop{3}\phi)(\delta, \delta, \delta) \\
&\le \frac{1}{6}\sup_{\param' \in \reals^{d}} \opnorm{\derivop{3}\phi(\param')}\twonorm{\delta}^{3}.
\]
Hence 
\[
\twonorm{\drift(\param) - \drift_{\mathrm{Laplace}}(\param)}^{2} 
&= \sum_{i=1}^{d} \bigg\{\frac{1}{2}\delta^{\top}(\hessian\partial_{i}\log\postdensity)(\map)\delta + R(\partial_{i} \log \postdensity, \param)\bigg\}^{2} \\
&= \sum_{i=1}^{d} \bigg\{\frac{1}{2}M_{i}^{\star}\twonorm{\delta}^{2} + \frac{1}{6}L_{i}\twonorm{\delta}^{3}\bigg\}^{2}.
\]
Let $\Lambda$ denote the diagonal matrix with $\Lambda_{ii} = \lambda_{i}~(i=1,\dots,d)$,
$X \dist \distNorm(0, \Lambda)$, and $\Param \dist \approxdist_{\mathrm{Laplace}}$.
Then $\twonorm{X}$ is equal in distribution to $\twonorm{\Param - \map}$. 
\begin{lemma} \label{lem:Gaussian-6-norm}
For $X$ defined above, $\EE(\twonorm{X}^{6}) = 8 \norm{\lambda}_{3}^{3} + 6\twonorm{\lambda}^{2}\onenorm{\lambda} + \onenorm{\lambda}^{3}$.
\end{lemma}
Using Lemmas \eqref{lem:Gaussian-norms} and \eqref{lem:Gaussian-6-norm}, we have
\[
\EE\left\{\twonorm{\drift(\Param) - \drift_{\mathrm{Laplace}}(\Param)}\right\}
&= \EE\left[\left\{\sum_{i=1}^{d} \bigg(\frac{1}{2}M_{i}^{\star}\twonorm{X}^{2} + \frac{1}{6}L_{i}\twonorm{X}^{3}\bigg)^{2}\right\}^{1/2}\right] \\
&\le \EE\left(\sum_{i=1}^{d} \frac{1}{2}M_{i}^{\star}\twonorm{X}^{2} + \frac{1}{6}L_{i}\twonorm{X}^{3}\right) \\
&\le \frac{1}{2}\onenorm{M^{\star}}C_{1}(\lambda) + \frac{1}{6}\onenorm{L}C_{2}(\lambda)^{3/2} \label{eq:p-1-Laplace-gradient-bound}
\]
and
\[
\EE\left\{\twonorm{\drift(\Param) - \drift_{\mathrm{Laplace}}(\Param)}^{2}\right\}
&= \EE\left\{\sum_{i=1}^{d} \bigg(\frac{1}{2}M_{i}^{\star}\twonorm{X}^{2} + \frac{1}{6}L_{i}\twonorm{X}^{3}\bigg)^{2}\right\} \\
&= \EE\left\{\sum_{i=1}^{d}\frac{1}{4}(M_{i}^{\star})^{2}\twonorm{X}^{4} + \frac{1}{6}M_{i}^{\star}L_{i}\twonorm{X}^{5} + \frac{1}{36}L_{i}^{2}\twonorm{X}^{6}\right\} \\
&\le \frac{1}{4}\twonorm{M^{\star}}^{2}C_{2}(\lambda)^{2} + \frac{1}{6}(M^{\star})^{\top}L C_{3}(\lambda)^{5/2} + \frac{1}{36}\twonorm{L}^{2}C_{3}(\lambda)^{3}. \label{eq:p-2-Laplace-gradient-bound}
\]
Using Eqs.~\eqref{eq:p-1-Laplace-gradient-bound} and \eqref{eq:p-2-Laplace-gradient-bound},
the result follows from \cref{thm:p-Wasserstein-sde-error} with $\mainmeas = \approxdist_{\mathrm{Laplace}}$
and $\estmeas = \postdist$.

\subsection*{Proof of \cref{lem:Gaussian-6-norm}}

Using $\EE(X_{i}^{6}) = 15\lambda_{i}^{3}$ and the identities from the proof of \cref{lem:Gaussian-norms},
we have
\[
\EE(\twonorm{X}^{6}) 
&= \EE\left\{\left(\textstyle\sum_{i=1}^{d}X_{i}^{2}\right)^{3}\right\} \\
&= \EE\left(\textstyle\sum_{i=1}^{d}X_{i}^{6} + 3\sum_{i=1}^{d}\sum_{j=1, j \ne i}^{d}X_{i}^{4}X_{j}^{2} + 6\sum_{i=1}^{d}\sum_{j=1}^{i-1}\sum_{k=1}^{j-1}X_{i}^{2}X_{j}^{2}X_{k}^{2}\right) \\
&= 15\textstyle\sum_{i=1}^{d}\lambda_{i}^{3} + 9\sum_{i=1}^{d}\sum_{j=1, j \ne i}^{d}\lambda_{i}^{2}\lambda_{j} + 6\sum_{i=1}^{d}\sum_{j=1}^{i-1}\sum_{k=1}^{j-1}\lambda_{i}\lambda_{j}\lambda_{k}  \\
&= 8\textstyle\sum_{i=1}^{d}\lambda_{i}^{3} + 6\sum_{i=1}^{d}\sum_{j=1}^{d}\lambda_{i}^{2}\lambda_{j} + \sum_{i=1}^{d}\sum_{j=1}^{d}\sum_{k=1}^{d}\lambda_{i}\lambda_{j}\lambda_{k}  \\
&= 8 \norm{\lambda}_{3}^{3} + 6\twonorm{\lambda}^{2}\onenorm{\lambda} + \onenorm{\lambda}^{3}.  \label{eq:sixth-moment}
\]
}

\bibliographystyle{biometrika}
\bibliography{references}

\bibliographystyle{numeric}
\bibliography{references}

\end{document}